\documentclass[review, sort&compress, 3p]{elsarticle}
\journal{Journal}

\usepackage{lineno, hyperref}
\modulolinenumbers[1]

\usepackage{makeidx}
\usepackage{color}

\usepackage{kpfonts}

\usepackage{amsmath}
\usepackage{amsfonts}
\usepackage{amssymb}
\usepackage{bm}
\usepackage{mathtools}

\usepackage{graphicx}
\usepackage{subfigure}
\usepackage{epsfig}
\usepackage{tikz}
\usetikzlibrary{snakes}
\usepackage{ifthen}
\usepackage{tikz-3dplot}

\usepackage{tabularx, multirow}	
\usepackage{footnote}	
\usepackage{array}

\usepackage{enumerate}
\usepackage{verbatim}
\usepackage{ulem}
\newcommand*{\abs}[1]{\left|#1\right|}

\newcommand{\diff}{\mathrm{d}}

\newcommand{\imag}{\mathrm{i}}

\newcommand{\trans}{\text{T}}
\newcommand{\herm}{\text{H}}

\newcounter{algleo}
\newlength{\lefttab}
\newlength{\numberoffset}
  {\trivlist
   \topsep=0pt\parsep=0pt\itemsep=0pt
   \addtolength{\lefttab}{1.25em}
   \addtolength{\numberoffset}{1.25em}
   \leftskip=\lefttab}%
  {\endtrivlist}

\begin{document}

\begin{frontmatter}
	\title{A nearly optimal explicitly-sparse representation 
	for oscillatory kernels with curvelet-like functions}
	\author{Yanchuang Cao}
	\author{Jun Liu}
	\author{Dawei Chen\corref{cor}}
	\ead{chen\_dawei@iapcm.ac.cn}

	\cortext[cor]{Corresponding author}
	\address{Institute of Applied Physics and Computational 
	Mathematics, Beijing 100094, China}

	\begin{abstract}
		A nearly optimal explicitly-sparse representation for oscillatory kernels
		is presented in this work by developing a curvelet based method.
		Multilevel curvelet-like functions are constructed as the transform of
		the original nodal basis. Then the system matrix in a new non-standard
		form is derived with respect to the curvelet basis, which would be nearly 
		optimally
		sparse due to the directional low rank property of the oscillatory kernel.
		Its sparsity is further enhanced via a-posteriori compression. Finally its
		nearly optimial log-linear computational complexity with controllable
		accuracy is demonstrated with numerical results.
		This explicitly-sparse representation is expected to lay ground to 
		future work related to fast direct solvers
		and effective preconditioners for high frequency problems. It may also be
		viewed as the generalization of wavelet based methods to high frequency
		cases, and used as a new wideband fast algorithm for wave problems.
	\end{abstract}
	\begin{keyword}
		Sparsification; curvelet; directional fast multipole method;
		high frequency; wavelet
	\end{keyword}
\end{frontmatter}


\section{Introduction}

In this work, $N$-body problems like
\begin{eqnarray} \label{eq_nbody}
	f_i = \sum_{j=1}^N K(\bm{x}_i, \bm{y}_j) \sigma_j,
	\quad \text{for }i = 1, 2, , , N
\end{eqnarray}
with the oscillatory kernel being
\begin{eqnarray} \label{eq_greenfunc}
	K(\bm{x}, \bm{y}) = \frac{\text{e}^{\imag \kappa r}}{4\pi r}, 
	\quad r = \abs{\bm{x} - \bm{y}}
\end{eqnarray}
or its derivatives are concerned.
Such computations may come from acoustic, electromagnetic, or elastodynamic
problems when solved with boundary integral equations or the 
Lippmann-Schwinger equation \cite{kirkup1998bem, caoyanchuang2015fdbem, 
wengchochew2013overview, jiangzhaoneng2019electromagnetic, 
caoyanchuang2019multidomain, chaillat2013elastodynamics}.

Two major difficulties may be encountered when solving large scale cases.
Firstly, the direct evaluation is prohibitive, 
since its computational complexity is $O(N^2)$ due to
its densely populated system matrix.
Secondly, iterative solvers which are generally used
often suffer from the bad conditioning of the system 
matrix, especially in cases with highly oscillatory kernels 
\cite{chaillat2013elastodynamics, liuyijun2019review, liufei2018preconditioner,
engquist2011sweeping}.
These difficulties can be overcome by constructing
approximate sparse representations for the system matrix and 
its inverse.

\emph{Data-sparse} representations for the system matrix are constructed
in many fast algorithms, that is, the system is represented in a low rank 
factorization form with only $O(N \log^\alpha N)_{\alpha \ge 0}$ elements.
For low frequency 
or even non-oscillatory cases which are mainly concerned in early research,
the kernel $K(\bm{x}, \bm{y})$ is asymptotically smooth when $\bm{x}$ 
lies far away from $\bm{y}$.
Therefore, $K(\bm{x}, \bm{y})$ can be approximated by a low-order expansion,
and low rank approximations for off-diagonal matrix blocks corresponding to 
far-field interactions can be constructed. 
Representative fast algorithms following this idea include the 
fast multipole method (FMM) and its variants \cite{rokhlin1985fmm, 
greengard1987fmm, yinglexing2004kifmm, fong2009bbfmm},
treecode algorithm \cite{barnes1986treecode},
$\mathcal{H}^2$-matrix \cite{hackbusch2002h2matrix},
adaptive cross approximation \cite{bebendorf2000aca},
and nested cross approximation \cite{bebendorf2012nca, gujjula2022nca}, etc.
One of the main differences between these algorithms is that they use
different functions in the low-order expansion of the kernel, or different 
algebraic schemes in deriving the low rank approximation for matrix 
blocks.
Each of them can achieve a data-sparse representation for the system matrix 
with only $O(N)$ elements, 
and bring down the computational complexity to linear for low
frequency cases.

For high frequency cases, however, data-sparse representations cannot be
obtained by these early fast algorithms, since
the approximate rank of the matrix block increases with its size
\cite{lucc1993fast}.
Thus other mathematical properties of the kernel should be considered
to develop fast algorithms. 
It is noticed that when the kernel is approximated with exponential
expansions, the translation matrices in FMM can be diagonalized 
even in high frequency, thus far field interactions can be computed 
efficiently. This leads to the diagonal form FMM 
\cite{rokhlin1993diagonal, greengard1997diagonalform} and
high frequency FMM \cite{darve2000fmm, nishimura2002fmm}.
To overcome its possible numerical instability in low frequency cases
\cite{coifman1993fmm, dembart1998fmm}, 
the wideband FMM is developed \cite{chenghongwei2006wbfmm, 
gumerov2009broadband}, which is essentially the hybrid of the 
classic and high frequency FMMs. That is, exponential 
expansions are only used when the oscillation over matrix blocks becomes
remarkable.
Consequently, the computational complexity is reduced to $O(N \log N)$
for high frequency cases, and remains $O(N)$ for low frequency cases.
Therefore, it provides a data-sparse representation
with low rank factorization and diagonalization for oscillatory kernels.

Another data-sparse representation with merely low rank factorization 
is proposed later based on
the \emph{directional low rank property} of the oscillatory kernel. That is, 
for $x \in \mathbb{X}$ and $y \in \mathbb{Y}$
satisfying the \emph{directional parabolic separation condition},
the kernel $K(\bm{x}, \bm{y})$ can be approximated by a low order expansion
\cite{engquist2007fda, messner2012dfmm}.
Therefore, the far-field interaction matrix would be always of approximately 
low rank if the far field is divided into multiple directional cones
when necessary. Their low rank approximations can be constructed
with the same manner as that in the low frequency fast algorithms.
In other words, low frequency fast algorithms can be generalized to 
high frequency cases by defining the directional cones, and data-sparse
representations of the system matrix can be achieved by carefully 
dividing the off-diagonal matrix blocks.
This idea is very attractive and leads to multiple wideband fast algorithms, 
including the directional FMM generalized from the kernel independent FMM
\cite{engquist2007fda, benson2014dfmm},
the directional FMM generalized from the black-box FMM
\cite{messner2012dfmm},
the dir$\mathcal{H}^2$-ACA method and 
the directional algebraic FMM based on nested cross approximation
\cite{bebendorf2015wideband, gujjula2022dafmm}. 
All of these directional algorithms can achieve the 
$O(N \log^\alpha N)_{\alpha \ge 0}$ computational complexity, which is 
comparable with the wideband FMM.

\emph{Explicitly-sparse} representations for the linear system are very useful
in the development of fast direct solvers and efficient preconditioners.
An representative one is the inverse FMM, which first transforms the data-sparse 
representation in FMM into an explicitly-sparse extended linear system, then 
compute its inverse efficiently \cite{ambikasaran2014ifmm, coulier2017ifmm}.
Unfortunately, it is based on classic FMM, thus
is only valid for low frequency cases.
For high frequency cases, a sparsify and sweep preconditioner is proposed recently 
for Lippmann-Schwinger equations, by which the computational cost can be reduced 
to nearly linear \cite{liufei2018preconditioner}.
But it requires computations on Cartesian grids over the entire computational
domain, thus its performance may be not satisfactory for problems with 
highly nonuniform distributed points, for example, boundary element analysis 
in which points are distributed only on the surface.

An explicitly-sparse system matrix can be achieved 
straightforwardly in the wavelet based method (WBM)
\cite{beylkin1991fwt, tausch2003multiscale, tausch2004wavelet}.
The boundary integral is discretized with wavelet-like functions which 
are compactly supported functions with ``vanishing'' or
``quasi-vanishing'' low order moments \cite{tausch2003multiscale, xiaojinyou2008wbem}.
Thus for low frequency problems, the interactions between well-separated 
wavelet-like functions are insignificant.
Only near-field interactions between wavelets have to be preserved in the system
matrix, making it explicitly-sparse with only $O(N)$ nonzero elements
\cite{beylkin1991fwt, tausch2004wavelet}. However, WBM is also only suitable
for low frequency cases \cite{huybrechs2004wavenumber, hawkins2007wbem, 
xiaojinyou2009wavelet}.

Both WBM and FMM are based on the fact that the kernel $K(\bm{x}, \bm{y})$ can be 
approximated by low order expansions for well separated $\bm{x}$ and $\bm{y}$.
Further study shows that they can transform into each other
\cite{xiaojinyou2009wavelet, xiaojinyou2010wavelet}, which means
it is possible to transform a data-sparse representation based on low rank
factorization into an explicitly-sparse representation with wavelet-like functions.
This finding evokes the idea that the data-sparse representation
in directional algorithms may also be transformed into an explicitly-sparse 
representation, and a new wideband fast algorithm with ``directional'' 
wavelet-like functions may be obtained.

Multiple directional wavelets have been proposed in the last two decades, among
them the \emph{curvelet} is a popular one and has gained great success in image
processing \cite{candes2000curvelet, candes2006fct, majianwei2010curvelet}.
Furthermore, mathematical study on its potential application in scientific 
computing shows that, the curvelet representation for oscillatory kernels
is optimally sparse, which consists only $O(N \log N)$ nonzero elements
\cite{candes2005curvelet}. Therefore, it is reasonable to expect that,
in the new wideband fast algorithm transformed from directional algorithms,
the new basis should be curvelet-like functions, thus the algorithm can be named
as \emph{curvelet based method} (CBM).

The aim of this paper is to develop a curvelet based method
providing a nearly optimal explicitly-sparse representation for 
oscillatory kernels. It is a transform of a directional FMM, and can also 
be viewed as a generalization of a WBM.
It is worth noting that in this work, we are not restricted to one 
particular low order expansion of the kernel function or an unique algorithm
constructing low rank approximations for matrix blocks,
thus it provides a framework transforming various directional algorithms
into CBMs.

This rest of the paper is organized as follows. In Section \ref{sec_wavelet}, 
a framework transforming FMM-like algorithms to WBM is discussed for low frequency 
problems. 
Then it is generalized to high frequency cases in Section \ref{sec_curvelet}, 
resulting in the CBM and explicitly-sparse representation for wideband problems.
Its sparsify is further enhanced by \emph{a-posteriori} compression in Section
\ref{sec_aposteriori}.
Finally its performance is studied numerically in Section \ref{sec_numres}.

\section{Wavelet compression method for low frequency cases}
\label{sec_wavelet}

Our CBM is transformed from a directional FMM. For low frequency cases, 
the directional FMM would degenerate into the corresponding low frequency 
FMM-like method, thus our CBM should degenerate into a WBM.
Therefore, before proposing CBM, we would like to provide
the framework transforming FMM-like algorithms to WBM for $N$-body problems
in this section.

\subsection{Basis and weight functions in $N$-body problems}

Generally the WBM is developed for Petrov-Galerkin
discretization of integrals \cite{tausch2003multiscale, tausch2004wavelet,
xiaojinyou2011wavelet}, in which the wavelet basis can be constructed by 
transforming the original nodal basis and weight functions.
For the $N$-body problem \eqref{eq_nbody}, it seems there are no basis or 
weight functions.
However, it can be rewritten into the following Petrov-Galerkin integral form
\begin{eqnarray} \label{eq_galerkin}
	f_i = \sum_{j=1}^N K(\bm{x}_i, \bm{y}_j) \sigma_j
	= \int_\Omega \delta_i(\bm{x}) \int_\Omega K(\bm{x}, \bm{y}) 
	\delta_j(\bm{y})\,\diff \bm{y} \,\diff \bm{x} \cdot \sigma_j, 
	\quad \bm{x_i} \in \Omega, \bm{y}_j \in \Omega,
\end{eqnarray}
where
\begin{eqnarray*}
	\delta_i(\bm{x}) = \delta(\bm{x}, \bm{x}_i), \quad 
	\delta_j(\bm{y}) = \delta(\bm{y}, \bm{y}_j).
\end{eqnarray*}
Thus the $N$-body problem can be viewed as integrals over $\Omega$ discretized
with Dirac-$\delta$ functions as basis and weight functions, and wavelets can be 
constructed as the linear combination of Dirac-$\delta$ functions.

\subsection{Wavelet compression for far-field interactions}
\label{subsec_sparsification}

For low frequency cases, the kernel can be approximated by low order
expansions 
\begin{eqnarray} \label{eq_loe}
	K(\bm{x}, \bm{y}) \approx \sum_r T_r(\bm{x}) \sum_s d_{rs} S_s(\bm{y}), 
	\quad \bm{x} \in \mathbb{X}, \bm{y} \in \mathbb{Y}
\end{eqnarray}
when $\mathbb{X}$ and $\mathbb{Y}$ are well separated
\begin{eqnarray}
	\text{dist}(\mathbb{X}, \mathbb{Y}) \ge 
	\eta \max(\text{diam} \mathbb{X}, \text{diam} \mathbb{Y}).
\end{eqnarray}
Therefore, far-field interactions in the Petrov-Galerkin disretization
can be evaluated by
\begin{eqnarray} \label{eq_slfmm}
	\begin{split}
		A_{ij} =& \int_\Omega w_i(\bm{x}) \int_\Omega K(\bm{x}, \bm{y}) 
		\chi_j(\bm{y})\,\diff \bm{y} \,\diff \bm{x} \\
		\approx& \sum_r \int_\Omega T_r(\bm{x}) w_i(\bm{x}) \,\diff \bm{x}
		\sum_s d_{rs} \int_\Omega S_s(\bm{y}) \chi_j(\bm{y}) \,\diff \bm{y},
		\quad \bm{x} \in \mathbb{X}, \bm{y} \in \mathbb{Y}.
	\end{split}
\end{eqnarray}
Denote 
\begin{eqnarray}
	M_{ri}(w) = \int_\Omega T_r(\bm{x}) w_i(\bm{x}) \,\diff \bm{x}, \quad
	M_{sj}(\chi) = \int_\Omega S_s(\bm{y}) \chi_j(\bm{y}) \,\diff \bm{y}
\end{eqnarray}
as moments of $w$ and $\chi$, respectively. 
Apparently $A_{ij} \approx 0$ if $\bm{M}(w_i) \approx 0$ or 
$\bm{M}(\chi_j) \approx 0$, i.e., when the moment of $w_i$ or $\chi_j$ 
approximately vanishes. Then the system matrix would be sparsified.

Wavelet-like functions are transformed from the original 
nodal functions $\{\chi\}$ and $\{w\}$. 
For simplicity they are also called wavelets in this paper.
For basis functions $\{ \chi \}$, 
evaluate the moment matrix $\bm{M}(\chi)$ and
calculate its singular value decomposition
\begin{eqnarray}
	\bm{M}(\chi) = \bm{U \Sigma} \bm{Q}^\herm,
\end{eqnarray}
in which $\bm{\Sigma} = \text{diag}(\sigma_0, \sigma_1, \cdots, \sigma_n)$ is 
the diagonal matrix consisting of singular values in the descending order.
Divide $\bm{\Sigma}$ and $\bm{Q}$ into columns corresponding to relatively large
singular values $\left\{ \sigma_i | \sigma_i \ge \varepsilon \sigma_0 \right\}$ and
ignorable ones $\left\{ \sigma_i | \sigma_i < \varepsilon \sigma_0 \right\}$, then
\begin{eqnarray} \label{eq_moment_split}
	\begin{split}
		\bm{M}(\chi) \left[ \bm{Q}_1, \bm{Q}_0 \right] 
		=& \bm{U} \left[ \bm{\Sigma}_1, \bm{\Sigma}_0 \right] 
		\approx \left[ \bm{U} \bm{\Sigma}_1, \bm{0} \right] \\
		:=& \left[ \bm{M}(\bm{\Phi}), \bm{M}(\bm{\Psi}) \right],
	\end{split}
\end{eqnarray}
with 
\begin{eqnarray}
	\left[ \bm{\Phi}, \bm{\Psi} \right] = 
	\left[ \bm{Q}_1^\trans \{\chi\}, \bm{Q}_0^\trans \{\chi\} \right].
\end{eqnarray}
That is, the original basis functions are transformed into a group of
scaling functions $\bm{\Phi} = \bm{Q}_1^\trans \{\chi\}$ without 
quasi-vanishing moments and
wavelets $\bm{\Psi} = \bm{Q}_0^\trans \{\chi\}$ with quasi-vanishing moments.

The number of scaling functions $\bm\Phi$ are limited and independent
with the matrix block's size. 
Assume the kernel is approximated by a $p$-term expansion, then the moment matrix
$\bm{M}(\chi)$ also consists of $p$ rows, thus the number of relatively
large singular values and number of scaling functions would not exceed $p$. Since
$p$ is independent of the number of points, there would be only $O(1)$ scaling
functions.

Wavelets for weight functions can be constructed in the similar manner.
Then the matrix block for far-field interactions \eqref{eq_slfmm} becomes
\begin{eqnarray}
	\begin{split}
		\bm{A} =& \left[\bm{M}(w)\right]^\trans \bm{D} \bm{M}(\chi)
		= \left\{\left[ \bm{M}(\bm{\Phi}_w),  \bm{M}(\bm{\Psi}_w) \right]
		\bm{Q}_w^\herm \right\}^\trans
		\bm{D} \left\{ \left[ \bm{M}({\bm{\Phi}_\chi}), \bm{M}({\bm{\Psi}_\chi}) \right] \bm{Q}_\chi^\herm \right\} \\
		\approx& \left[ \bm{M}(\bm{\Phi}_w) \bm{Q}_{w, 1}^\herm \right]^\trans
		\bm{D} \left[ \bm{M}({\bm{\Phi}_\chi}) \bm{Q}_{\chi, 1}^\herm \right]
		= \overline{\bm{Q}_{w, 1}} \left\{ \left[\bm{M}(\bm{\Phi}_w)\right]^\trans
		\bm{D} \bm{M}(\bm{\Phi}_\chi) \right\} \bm{Q}_{\chi, 1}^\herm \\
		:=& \overline{\bm{Q}_{w, 1}} \tilde{\bm{A}} \bm{Q}_{\chi, 1}^\herm
	\end{split}
\end{eqnarray}
with $\bm{D} = \left[ d_{rs} \right]$.
Hence the matrix block $\bm{A}$ can be represented with
\begin{eqnarray}
	\tilde{\bm{A}} = \bm{Q}_{w, 1}^\trans \bm{A} \bm{Q}_{\chi, 1}
	= \left[\bm{M}(\bm{\Phi}_w)\right]^\trans \bm{D} \bm{M}(\bm{\Phi}_\chi)
\end{eqnarray}
which only consists of interactions of scaling functions
$\bm{\Phi}_w$ and $\bm{\Phi}_\chi$, and there are only $O(1)$ nonzero elements.

Notice Eq. \eqref{eq_slfmm} is actually the translation in single level FMM, in 
which $\bm{M}(\chi), \bm{D}, \left[\bm{M}(w)\right]^\trans$ are called the 
source-to-moment (S2M), moment-to-local (M2L), and local-to-target (L2T)
translation matrices, respectively.
That is, when transforming FMM into WBM, the moment matrix of basis functions
can be defined the same as the S2M matrix in FMM, while the 
moment matrix for weight functions should be defined as the transpose of
the L2T matrix.

\subsection{Multilevel wavelet construction}

An explicitly-sparse representation for the system matrix with $O(N)$ nonzero
elements can be achieved with multilevel wavelet basis. 
They can be constructed on a balanced octree, which is commonly used in WBMs.
First define a Level-0 cube containing all points. Then for currently the finest 
level, if there are cubes containing more than $N_p$ points, each cube
in this level would be divided into eight child cubes in the next level.
The subdivision is continued until the number of nodes in each leaf cube does
not exceed the predetermined number $N_p$.
For each cube $C$, define its near field $N^C$ as the union of cubes adjacent
with $C$. The rest are defined as its far field $F^C$, and the interaction field 
is defined as $I^C = N^P \backslash N^C$ with $P$ being the parent of $C$.

The construction of multilevel wavelets starts from the finest $L$-th level.
For each leaf cube $C$, compute the moment matrix for nodal basis inside $C$.
Then transform the nodal basis into wavelets $\bm{\Psi}_L$ and scaling functions 
$\bm{\Phi}_L$, and compute the moments of scaling functions $\bm{M}_L(\bm{\Phi}_L)$
in the $L$-th level.

For each cube $C$ in the $l$-th level with $l<L$, the scaling functions in its 
child cubes are collected and transformed into wavelets and scaling functions 
in the $l$-th level. First translate the moments of scaling functions in its 
child cubes $\bm{M}_{l+1}(\bm{\Phi}_{l+1})$ into moments on the $l$-th level 
by the moment-to-moment (M2M) translation
\begin{eqnarray}
	\bm{M}_l(\bm{\Phi}_{l+1}) = \mathcal{M}_l \bm{M}_{l+1}(\bm{\Phi}_{l+1}),
\end{eqnarray}
then compute its singular value decomposition to get
wavelets $\bm{\Psi}_l$ and scaling functions $\bm{\Phi}_l$ 
in the $l$-th level and their moments $\bm{M}_l(\bm{\Psi}_l)$ and 
$\bm{M}_l(\bm{\Phi}_l)$. This is done recursively until the $h$-th level
is reached on which all cubes lie in the near field or interaction field of 
each other.

Now let's discuss the multilevel wavelet construction for weight functions.
Notice that the translations in multilevel FMM can be expressed as
\begin{eqnarray}
	\bm{b} = \mathcal{TLDMS} \bm{x},
\end{eqnarray}
with $\mathcal{S}, \mathcal{M}, \mathcal{D}, \mathcal{L}, \mathcal{T}$
denote the S2M, M2M, M2L, L2L, L2T translations, respectively.
Comparing with the single level FMM $\bm{b} = \mathcal{TDS}\bm{x}$,
$\mathcal{MS}$ actually takes the role of
$\left[ \bm{M}_{L-1}(\bm{\Phi}_\chi) \right]$ 
at the $(L-1)$-th levels with $\mathcal{S} = \bm{M}_L (\chi)$, 
and $\mathcal{TL}$ takes the role of
$\left[ \bm{M}_{L-1}(\bm{\Phi}_w) \right]^\trans$ 
with $\mathcal{T} = \left[ \bm{M}_L(w)\right]^\trans$.
Therefore, for weight functions,
\begin{eqnarray}
	\bm{M}_l(\bm{\Phi}_{l+1}) = \mathcal{L}_l^\trans \bm{M}_{l+1}(\bm{\Phi}_{l+1}).
\end{eqnarray}
That is, in the multilevel wavelet construction for weight functions, 
the moments of scaling functions should be translated
by the transpose of the L2L translation matrix in FMM.

The construction of multilevel wavelet basis from
original nodal basis $\bm{X} = \{\chi\}$ can be depicted as
\begin{eqnarray*}
	\begin{matrix}
		\bm{X} & \rightarrow & \bm{\Phi}_L & \rightarrow &
		\bm{\Phi}_{L-1}
		& \rightarrow & \cdots & \rightarrow & \bm{\Phi}_h \\
		~ & \searrow & ~ & \searrow & ~ & \searrow & ~ & \searrow \\
		~ & ~ & \bm{\Psi}_L & ~ & \bm{\Psi}_{L-1} & ~ & \cdots & ~ &
		\bm{\Psi}_h
	\end{matrix}
\end{eqnarray*}

\subsection{Linear system with respect to the wavelet basis}

There are two kinds of system matrix representations with respect to the multilevel wavelet
basis, namely the standard form and the non-standard form. The non-standard form
is used in this work, in which there are only interactions of functions
on the same level, thus it is much simpler. Moreover, it is more efficient
than the standard form \cite{beylkin1991fwt}.

The system matrix can be transformed into wavelet representations
with wavelets in the $L$-th level
\begin{eqnarray}
	\bm{A} = \left[ \begin{matrix}
		\overline{\bm{Q}_{w,0,L}} & \overline{\bm{Q}_{w,1,L}}
	\end{matrix} \right]
	\left[ \begin{matrix}
		\bm{A}_L^{\psi, \psi} & \bm{A}_L^{\psi, \phi} \\
		\bm{A}_L^{\phi, \psi} & \bm{A}_L^{\phi, \phi}
	\end{matrix} \right]
	\left[ \begin{matrix}
		\bm{Q}_{\chi,0,L}^\herm \\
		\bm{Q}_{\chi,1,L}^\herm
	\end{matrix} \right].
\end{eqnarray}
Far-field interactions with wavelets are always insignificant due to their 
quasi-vanishing moments, thus $\bm{A}_L^{\psi, \psi}, \bm{A}_L^{\psi, \phi}, 
\bm{A}_L^{\phi, \psi}$ are sparse, and only $\bm{A}_L^{\phi, \phi}$ is fully 
populated. $\bm{A}_L^{\phi, \phi}$ can be further sparsified with wavelet basis in 
the $(L-1)$-th level, which gives
\begin{eqnarray}
	\begin{split}
		\bm{A} =& \left[ \begin{matrix}
			\overline{\bm{Q}_{w,0,L}} & \overline{\bm{Q}_{w,1,L}}
		\end{matrix} \right]
		\left[ \begin{matrix}
			\bm{A}_L^{\psi, \psi} & \bm{A}_L^{\psi, \phi} \\
			\bm{A}_L^{\phi, \psi} & \bm{A}_L^{\phi, \phi}
		\end{matrix} \right]
		\left[ \begin{matrix}
			\bm{Q}_{\chi,0,L}^\herm \\
			\bm{Q}_{\chi,1,L}^\herm
		\end{matrix} \right] \\
		=& \left[ \begin{matrix}
			\overline{\bm{Q}_{w,0,L}} & \overline{\bm{Q}_{w,1,L}}
		\end{matrix} \right]
		\left[ \begin{matrix}
			\bm{A}_L^{\psi, \psi} & \bm{A}_L^{\psi, \phi} \\
			\bm{A}_L^{\phi, \psi} &
			\left[ \begin{matrix}
				\overline{\bm{Q}_{w,0,L-1}} & \overline{\bm{Q}_{w,1,L-1}}
			\end{matrix} \right]
			\left[ \begin{matrix}
				\bm{A}_{L-1}^{\psi, \psi} & \bm{A}_{L-1}^{\psi, \phi} \\
				\bm{A}_{L-1}^{\phi, \psi} & \bm{A}_{L-1}^{\phi, \phi}
			\end{matrix} \right]
			\left[ \begin{matrix}
				\bm{Q}_{\chi,0,L-1}^\herm \\
				\bm{Q}_{\chi,1,L-1}^\herm
			\end{matrix} \right]
		\end{matrix} \right]
		\left[ \begin{matrix}
			\bm{Q}_{\chi,0,L}^\herm \\
			\bm{Q}_{\chi,1,L}^\herm
		\end{matrix} \right] \\
		=& \left[ \begin{matrix}
			\overline{\bm{Q}_{w,0,L}} & \overline{\bm{Q}_{w,1,L}} & \overline{\bm{Q}_{w,1,L}}\, \overline{\bm{Q}_{w,0,L-1}} & \overline{\bm{Q}_{w,1,L}}\, \overline{\bm{Q}_{w,1,L-1}} 
		\end{matrix} \right]
		\left[ \begin{matrix}
			\bm{A}_L^{\psi, \psi} & \bm{A}_L^{\psi, \phi} & ~ & ~ \\
			\bm{A}_L^{\phi, \psi} & ~ & ~ & ~ \\
			~ & ~ & \bm{A}_{L-1}^{\psi, \psi} & \bm{A}_{L-1}^{\psi, \phi} \\
			~ & ~ & \bm{A}_{L-1}^{\phi, \psi} & \bm{A}_{L-1}^{\phi, \phi}
		\end{matrix} \right]
		\left[ \begin{matrix}
			\bm{Q}_{\chi,0,L}^\herm \\
			\bm{Q}_{\chi,1,L}^\herm \\
			\bm{Q}_{\chi,0,L-1}^\herm \bm{Q}_{\chi,1,L}^\herm \\
			\bm{Q}_{\chi,1,L-1}^\herm \bm{Q}_{\chi,1,L}^\herm
		\end{matrix} \right].
	\end{split}
\end{eqnarray}
This is done recursively until the highest level. 
Consequently, the matrix-vector multiplication $\bm{b = Ax}$ is transformed 
into 
\begin{eqnarray}
	\bm{b} = \overline{\bm{Q}_w} \tilde{\bm{A}}_\text{ns} \bm{Q}_\chi^\herm \bm{x}.
\end{eqnarray}
The system matrix in the non-standard form
\begin{eqnarray} \label{eq_ans_wavelet}
	\tilde{\bm{A}}_\text{ns} = \left[ \begin{matrix}
		\bm{A}_L^{\psi, \psi} & \bm{A}_L^{\psi, \phi} & ~ & ~ & ~ & ~ & ~\\
		\bm{A}_L^{\phi, \psi} & ~ & ~ & ~ & ~ & ~ & ~\\
		~ & ~ & \ddots & ~ & ~ & ~ & ~\\
		~ & ~ & ~ & \bm{A}_{h+1}^{\psi, \psi} & \bm{A}_{h+1}^{\psi, \phi} & ~ & ~ \\
		~ & ~ & ~ & \bm{A}_{h+1}^{\phi, \psi} & ~ & ~ & ~ \\
		~ & ~ & ~ & ~ & ~ & \bm{A}_h^{\psi, \psi} & \bm{A}_h^{\psi, \phi} \\
		~ & ~ & ~ & ~ & ~ & \bm{A}_h^{\phi, \psi} & \bm{A}_h^{\phi, \phi}
	\end{matrix} \right]
\end{eqnarray}
is optimally sparse, since far-field interactions in
$\bm{A}^{\psi, \psi}, \bm{A}^{\psi, \phi}, \bm{A}^{\phi, \psi}$
can be discarded and the size of $\bm{A}_h^{\phi, \phi}$ is limited.
\begin{eqnarray}
	\begin{split}
		\bm{Q}_\ast = &\left[ 
		\bm{Q}_{\ast,0,L},\quad
		\bm{Q}_{\ast,1,L},\quad
		\bm{Q}_{\ast,1,L}\bm{Q}_{\ast,0,L-1},\quad
		\bm{Q}_{\ast,1,L}\bm{Q}_{\ast,1,L-1},\quad
		\cdots,\quad \right. \\
		& \left. \bm{Q}_{\ast,1,L}\bm{Q}_{\ast,1,L-1}\cdots\bm{Q}_{\ast,1,h+1}\bm{Q}_{\ast,0,h},\quad
		\bm{Q}_{\ast,1,L}\bm{Q}_{\ast,1,L-1}\cdots\bm{Q}_{\ast,1,h+1}\bm{Q}_{\ast,1,h}
		\right] \\
		= & \left[
		\bm{Q}_{\ast,0,L},\quad
		\bm{Q}_{\ast,1,L},\quad
		\bm{Q}_{\ast,1,L}\bm{Q}_{\ast,0,L-1},\quad
		\bm{Q}_{\ast,1,L}\bm{Q}_{\ast,1,L-1},\quad
		\cdots, \quad
		\prod_{l=L}^{h+1} \bm{Q}_{\ast,1,l} \cdot \bm{Q}_{\ast,0,h},\quad
		\prod_{l=L}^{h+1} \bm{Q}_{\ast,1,l} \cdot \bm{Q}_{\ast,1,h}
		\right]
	\end{split}
\end{eqnarray}
is the transform matrix between the original nodal basis $\bm{X}$ and wavelet
basis
\begin{eqnarray}
	\bm{\Psi}_\ast = \left\{ \bm{\Phi}_{\ast,L}, \bm{\Psi}_{\ast,L}, 
	\bm{\Phi}_{\ast,L-1}, \bm{\Psi}_{\ast,L-1}, \cdots, 
	\bm{\Phi}_{\ast,h}, \bm{\Psi}_{\ast,h}
	\right\},
\end{eqnarray}
with $\ast$ in the subscript denoting $w$ for weight functions or $\chi$ for basis 
functions.
Then the matrix vector 
multiplication $\bm{b = Ax}$ can be computed with linear computational cost
by the following three steps in WBM:
\begin{enumerate}
	\item Compute the coefficients $\tilde{\bm{x}}$ of wavelet basis functions
		$\bm{\Psi}_\chi$ with fast wavelet transform
		$\tilde{\bm{x}} = \bm{Q}_\chi^\herm \bm{x}$;
	\item Evaluate coefficients $\tilde{\bm{b}}$ of wavelet weight functions
		$\bm{\Psi}_w$ by $\tilde{\bm{b}} = \tilde{\bm{A}}_\text{ns} \tilde{\bm{x}}$;
	\item Compute the target values with inverse fast wavelet transform
		$\bm{b} = \overline{\bm{Q}_w} \tilde{\bm{b}}$.
\end{enumerate}

The non-standard system matrix can be computed in a upward pass.
Matrix blocks corresponding to a cube $C$ and $B \in N^C$
is computed by
\begin{eqnarray} \label{eq_mateval}
    \left[ \begin{matrix}
		\bm{A}_{B,C}^{\phi, \phi} & \bm{A}_{B,C}^{\phi, \psi}  \\
		\bm{A}_{B,C}^{\psi, \phi} & \bm{A}_{B,C}^{\psi, \psi}
	\end{matrix} \right] 
	= \bm{Q}_{L,w}^\trans 
	\bm{A}_{B,C}
	\bm{Q}_{L,\chi}.
\end{eqnarray}
If $C$ and $B$ lie in the finest level, $\bm{A}_{B,C}$ consists of interactions
of the original nodal basis and weight functions, thus it is actually the 
source-to-target (S2T) translation matrix in FMM.
If $C$ and $B$ lie in a higher $l$-th level,
$\bm{A}_{B,C}$ consists of interactions of all the $\phi$-interactions of
their children, which includes near field and interaction field interactions
in the $(l+1)$-th level. The near field interactions have already been computed 
by the transformation \eqref{eq_mateval} in the $(l+1)$-th level. 
The interaction field interactions can be computed efficiently with the moment
matrix $\bm{M}(\bm{\Phi})$ and the M2L translation
\begin{eqnarray} \label{eq_phi_interaction}
	\bm{A}^{\phi,\phi} = 
	\bm{Q}_w^\trans \bm{A} \bm{Q}_\chi
	= \bm{Q}_w^\trans \left[ \bm{M}(w) \right]^\trans \bm{D}
	\bm{M}(\chi) \bm{Q}_\chi
	= \left[ \bm{M}(\bm{\Phi}_w) \right]^\trans \bm{D} \bm{M}(\bm{\Phi}_\chi).
\end{eqnarray}
The matrix computation only requires computations in near field and interaction
field, thus its computational cost is linear.


\section{Curvelet compression for high frequency cases}
\label{sec_curvelet}

Now we transform the directional FMM to a CBM.
Our CBM also consists of computations in low and high frequency regime, 
and the computations in the low frequency regime is the same with the WBM.
Thus in this section, we mainly discuss the computations in the high frequency 
regime.

\subsection{Multilevel curvelet construction}

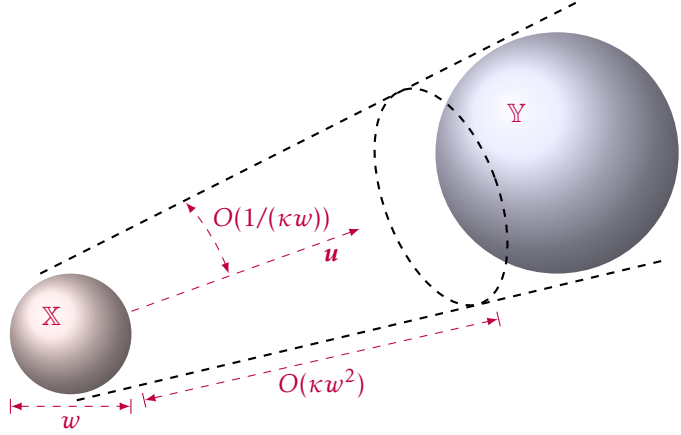
\begin{figure}[ht]
	\centering
	\begin{tikzpicture}[>=latex, scale=0.8]
		\shade[ball color=red!10] (0,0) circle(1cm) node[purple, above left]{$\mathbb{X}$};
		\shade[ball color=blue!10] (8,3) circle(2cm) node[purple, above left=8pt]{$\mathbb{Y}$};
		\draw[thick, dashed] (-0.5, 1) -- +(27:10);
		\draw[thick, dashed] (0.1, -1.1) -- +(13.5:10);
		\draw[thick, dashed, rotate=20] (6.5,0.05) ellipse (0.94 and 1.88);
		\draw[purple, dashed, ->] (1, 0.375) --  node [very near end, below]{$\bm{u}$} +(20:4);
		\draw[purple, dashed, <->] (2.6, 1.0) arc (20:42:3.6);
		\draw[purple] (3.3, 1.9) node {$O(1/(\kappa w))$};
		\draw[purple, dashed, |<->|] (-1, -1.2) -- node[below]{$w$} (1, -1.2);
		\draw[purple, dashed, |<->|] (1.2, -1.1) -- node[below]{$O(\kappa w^2)$} +(13:6.0);
	\end{tikzpicture}
	\caption{Parabolic separation condition.}
	\label{fig_parasep}
\end{figure}

Our CBM is also based on the directional low rank property of the oscillatory 
kernel. That is, the kernel can be approximated by low order expansion 
\begin{eqnarray*}
	K(\bm{x}, \bm{y}) \approx \sum_r T_r(\bm{x}) \sum_s d_{rs} S_s(\bm{y}), 
	\quad \bm{x} \in \mathbb{X}, \bm{y} \in \mathbb{Y}
\end{eqnarray*}
when $\mathbb{X}$ and $\mathbb{Y}$ satisfies the parabolic separation condition,
as illustrated in Fig. \ref{fig_parasep}.
According to the theory in Section \ref{subsec_sparsification}, 
far-field interactions between $\mathbb{X}$ and $\mathbb{Y}$ can be
sparsified by using basis and weight functions with quasi-vanishing moments.
Notice that the expansion is only valid for a directional cone, and 
$T_r(\bm{x})$ and $S_s(\bm{y})$ differ from cone to cone, thus
the moments and wavelet-like functions are directional.
In the following, these directional wavelet-like functions are called 
curvelet-like functions, or curvelets for simplicity in this work.

In the multilevel curvelet construction, the octree is divided into low and 
high frequency regimes. The low frequency 
regime consists of low level cubes whose size is less than the wavelength,
and the rest lie in the high frequency regime.
Denote the highest level in the low frequency regime is the $h_l$-th level.
In order to ensure $h_l < L$ so that the low frequency regime is not empty, 
cubes larger than the wavelength are always subdivided during the octree 
construction.

For a cube $C$, if it lies in the low frequency regime, its near field, far field,
and interaction field are defined the same as that in WBM. If it lies in 
the high frequency regime, its near field $N^C$ is defined as the union of cubes
$B$ whose distance does not exceed $O(kw^2)$, and the rest is defined as its
far field $F^C$. Its interaction field $I^C = N^P \backslash N^C$ is divided into
multiple directional cones with spanning angles being $O(1/(\kappa w))$.
Thus each directional cone is parabolic separated from $C$.

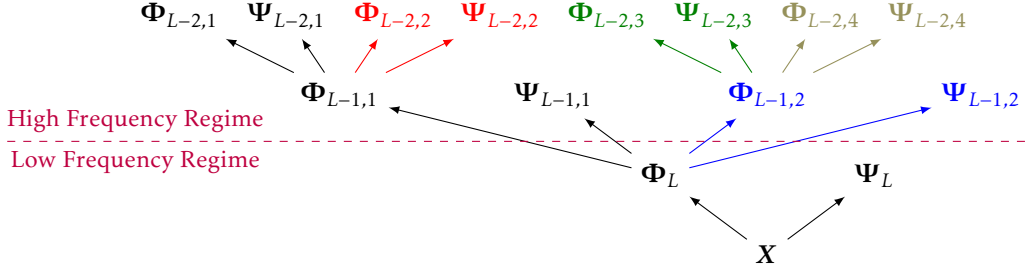
\begin{figure}[ht]
	\centering
	\begin{tikzpicture}[>=latex, level distance=30pt]
		\tikzstyle{level 1} = [sibling distance = 80pt]
		\tikzstyle{level 2} = [sibling distance = 80pt]
		\tikzstyle{level 3} = [sibling distance = 40pt]
		\node {$\bm{X}$} [grow'=up]
		  child {node{$\bm{\Phi}_L$} edge from parent[->]
			child {node{$\bm{\Phi}_{L-1,1}$} edge from parent[->]
			  child {node{$\bm{\Phi}_{L-2,1}$} edge from parent[->]}
			  child {node{$\bm{\Psi}_{L-2,1}$} edge from parent[->]}
			  child[red] {node{$\bm{\Phi}_{L-2,2}$} edge from parent[->]}
			  child[red] {node{$\bm{\Psi}_{L-2,2}$} edge from parent[->]}
			}
			child {node{$\bm{\Psi}_{L-1,1}$} edge from parent[->]}
			child[blue] {node{$\bm{\Phi}_{L-1,2}$} edge from parent[->]
			  child[green!50!black] {node{$\bm{\Phi}_{L-2,3}$} edge from parent[->]}
			  child[green!50!black] {node{$\bm{\Psi}_{L-2,3}$} edge from parent[->]}
			  child[yellow!50!black] {node{$\bm{\Phi}_{L-2,4}$} edge from parent[->]}
			  child[yellow!50!black] {node{$\bm{\Psi}_{L-2,4}$} edge from parent[->]}
			}
			child[blue] {node{$\bm{\Psi}_{L-1,2}$} edge from parent[->]}
		  }
		  child {node{$\bm{\Psi}_L$} edge from parent[->]};
		  \draw[dashed, purple] (-10,1.5) -- node[very near start, above] {\small High Frequency Regime} node[very near start, below] {\small Low Frequency Regime} (3.5,1.5);
	\end{tikzpicture}
	\caption{Construction of multilevel curvelets, where subscripts $(l, \gamma)$ 
	denotes the $\gamma$-th directional cone in the $l$-th level.}
	\label{fig_construction}
\end{figure}

The construction of multilevel curvelets from original nodal basis $\bm{X}$ 
are depicted in Fig. \ref{fig_construction}.
The original nodal basis are firstly transformed into wavelet basis
in the low frequency regime, then the scaling functions are further transformed
into curvelets in the high frequency regime.
For each directional cone, a group of curvelets and directional scaling
functions are constructed.
For the $(h_l-1)$-th level, i.e., the finest level in the high frequency regime, 
curvelets in the $\bm{u}$-th directional cone are transformed from non-directional 
scaling functions in the $h_l$-th level.
For a higher $l$-th level in the high frequency regime with $l < h_l$, curvelets 
in the $\bm{u}_l$-th directional cone are transformed from directional scaling 
functions in the $\bm{u}_{l+1}$-th directional cone in the $(l+1)$-th level, 
with $\bm{u}_l$ enclosed in the $\bm{u}_{l+1}$-th directional cone.

Notice there are $O(N \log N)$ directional cones for all cubes in the high frequency 
regime \cite{engquist2007fda}. For each directional cone, the curvelet basis 
is constructed via M2M and SVD of the moment matrix. Thus it is easy to find that
the computational complexity of the multilevel curvelet construction is $O(N \log N)$.


\subsection{Linear system in the new non-standard form}

Now let's construct the explicitly-sparse representation of the system matrix
with respect to the curvelet basis.
The matrix sparsification in the low frequency regime is the same with that in the 
WBM, resulting in the transformed system matrix like Eq. \eqref{eq_ans_wavelet}.
Then there would be only one densely populated matrix block 
$\bm{A}_{h_l}^{\phi, \phi}$
corresponding to interactions of all scaling functions at the $h_l$-th level,
and its size grows with the wavenumber $\kappa$.
Thus the key of curvelet sparsification is to further sparsify
$\bm{A}_{h_l}^{\phi, \phi}$
with curvelet basis in the high frequency regime.

A straightforward idea to sparsify $\bm{A}_{h_l}^{\phi, \phi}$ may be to apply
the sparsification scheme in WBM straightforwardly, i.e., to transform
it into near-field interactions of curvelets and directional scaling functions.
However, there are multiple groups of curvelet basis, each for a directional cone.
$\bm{A}_{h_l}^{\phi, \phi}$ can be transformed by
\begin{eqnarray}
	\bm{A}_{h_l}^{\phi, \phi} = 
	\left[ \begin{matrix}
		\overline{\bm{Q}_{w,0,h_l-1,\gamma}} &
		\overline{\bm{Q}_{w,1,h_l-1,\gamma}}
	\end{matrix} \right]
	\left[ \begin{matrix}
		\bm{A}_{h_l-1,\gamma}^{\psi, \psi} & \bm{A}_{h_l-1,\gamma}^{\psi, \phi} \\
		\bm{A}_{h_l-1,\gamma}^{\phi, \psi} & \bm{A}_{h_l-1,\gamma}^{\phi, \phi}
	\end{matrix} \right]
	\left[ \begin{matrix}
		\bm{Q}_{\chi,0,h_l-1,\gamma}^\herm \\
		\bm{Q}_{\chi,1,h_l-1,\gamma}^\herm
	\end{matrix} \right]
\end{eqnarray}
for each $\gamma$-th group. We cannot determine which group should be used.
More importantly, the transformed matrix blocks
$\bm{A}_{h_l,\gamma}^{\psi, \psi}, \bm{A}_{h_l,\gamma}^{\psi, \phi},
\bm{A}_{h_l,\gamma}^{\phi, \psi}$
are not sparse any more, since only $\psi$-interactions 
in the $\gamma$-th directional cone is ignorable, while that in other cones
are still significant. Therefore, we must find another scheme
to sparsify $\bm{A}_{h_l}^{\phi, \phi}$ and construct a new 
non-standard form to achieve the sparse representation for the system matrix.

In this work, the far-field interactions in $\bm{A}_{h_l}^{\phi, \phi}$ are 
sparsified cone-by-cone. Notice that $\bm{A}_{h_l}^{\phi, \phi}$ can be divided
into the near-field and interaction-field interactions at the $h_l$-th level
$\bm{A}_{h_l,\text{NI}(h_l)}^{\phi, \phi}$
and the interaction-field interactions on higher levels,
and the later is further divided into interactions with multiple directional cones.
That is,
\begin{eqnarray} \label{eq_cvdiv}
	\bm{A}_{h_l}^{\phi, \phi}
	= \bm{A}_{h_l,\text{NI}(h_l)}^{\phi, \phi} + 
	\sum_{l=h_l-1}^h \sum_\gamma^{\Gamma_l} \bm{A}_{h_l, \text{I}(l, \gamma)}^{\phi, \phi},
\end{eqnarray}
where $\text{I}(l,\gamma)$ denotes the $\gamma$-th directional cone in the 
interaction field in the $l$-th level, $\Gamma_l$ is the number of cones in the
$l$-th level, and
\begin{eqnarray*}
	\bm{A}_{h_l, \text{I}(l, \gamma)}^{\phi, \phi} = 
	\int_\mathbb{X} \bm{\Phi}_{w, h_l}(\bm{x}) \int_\mathbb{Y} K(\bm{x}, \bm{y}) 
	\bm{\Phi}_{\chi, h_l}(\bm{y}) \,\diff\bm{y} \,\diff\bm{x}
\end{eqnarray*}
with $\mathbb{X}$ being the cubes in the $l$-th level, and $\mathbb{Y}$ being
their $\gamma$-th directional cones.

$\bm{A}_{h_l, \text{I}(l, \gamma)}^{\phi, \phi}$ can be approximated by 
low order expansion since the kernel is of directional low rank.
First transform the $\phi_{h_l}$-interactions into $\phi_{h_l-1}$-interactions
\begin{eqnarray} \label{eq_cvcmp}
	\bm{A}_{h_l, \text{I}(l, \gamma)}^{\phi, \phi} = 
	\left[ \begin{matrix}
		\overline{\bm{Q}_{w,0,h_l-1,\gamma'}} &
		\overline{\bm{Q}_{w,1,h_l-1,\gamma'}}
	\end{matrix} \right]
	\left[ \begin{matrix}
		\bm{A}_{h_l-1,\gamma'}^{\psi, \psi} & \bm{A}_{h_l-1,\gamma'}^{\psi, \phi} \\
		\bm{A}_{h_l-1,\gamma'}^{\phi, \psi} & \bm{A}_{h_l-1,\gamma'}^{\phi, \phi}
	\end{matrix} \right]
	\left[ \begin{matrix}
		\bm{Q}_{\chi,0,h_l-1,\gamma'}^\herm \\
		\bm{Q}_{\chi,1,h_l-1,\gamma'}^\herm
	\end{matrix} \right]
\end{eqnarray}
where the subscript $(h_l-1, \gamma')$ denotes the directional cone in $(h_l-1)$-th 
level enclosing $\text{I}(l,\gamma)$.
Matrix blocks
$\bm{A}_{h_l-1,\gamma'}^{\psi, \psi}$, $\bm{A}_{h_l-1,\gamma'}^{\psi, \phi}$ and
$\bm{A}_{h_l-1,\gamma'}^{\phi, \psi}$
are ignorable due to the quasi-vanishing directional moments of curvelets, thus
Eq. \eqref{eq_cvcmp} becomes
\begin{eqnarray}
	\bm{A}_{h_l, \text{I}(l, \gamma)}^{\phi, \phi} = 
	\overline{\bm{Q}_{w,1,h_l-1,\gamma'}}
	\bm{A}_{h_l-1,\gamma'}^{\phi, \phi}
	\bm{Q}_{\chi,1,h_l-1,\gamma'}^\herm.
\end{eqnarray}
That is, $\phi_{h_l}$-interactions are transformed into $\phi_{h_l-1}$-interactions,
and the matrix block is compressed.
This is done recursively until the $l$-th level, which gives
\begin{eqnarray} \label{eq_cvcmprec}
	\begin{split}
		\bm{A}_{h_l, \text{I}(l, \gamma)}^{\phi, \phi} = &
		\left(
		\overline{\bm{Q}_{w,1,h_l-1,\gamma'}}\,
		\overline{\bm{Q}_{w,1,h_l-2,\gamma'}}
		\cdots
		\overline{\bm{Q}_{w,1,l,\gamma}}
		\right)
		\bm{A}_{l,\gamma}^{\phi, \phi}
		\left(
		\bm{Q}_{\chi,1,l,\gamma}^\herm
		\cdots
		\bm{Q}_{\chi,1,h_l-2,\gamma'}^\herm
		\bm{Q}_{\chi,1,h_l-1,\gamma'}^\herm
		\right) \\
		= & \prod_{\lambda=h_l-1}^l \overline{\bm{Q}_{w, 1, \lambda, \gamma'}} \cdot
		\bm{A}_{l,\gamma}^{\phi, \phi}
		\cdot \prod_{\lambda=l}^{h_l-1} \bm{Q}_{\chi, 1, \lambda, \gamma'}^\herm.
	\end{split}
\end{eqnarray}

Substitute Eq. \eqref{eq_cvcmprec} into \eqref{eq_cvdiv}, we get
\begin{eqnarray}
	\bm{A}_{h_l}^{\phi, \phi}
	= \bm{A}_{h_l,\text{NI}(h_l)}^{\phi, \phi} + 
	\sum_{l=h_l-1}^h \sum_\gamma^{\Gamma_l} 
	\prod_{\lambda=h_l-1}^l \overline{\bm{Q}_{w, 1, \lambda, \gamma'}} \cdot
	\bm{A}_{l,\gamma}^{\phi, \phi}
	\cdot \prod_{\lambda=l}^{h_l-1} \bm{Q}_{\chi, 1, \lambda, \gamma'}^\herm.
\end{eqnarray}
It can be written into the form
\begin{eqnarray} \label{eq_ahlphiphi}
	\bm{A}_{h_l}^{\phi, \phi} = \overline{\bm{Q}_{w, \text{HFR}}} \cdot
	\bm{A}_{\text{ns},\text{HFR}} \cdot \bm{Q}_{\chi, \text{HFR}}^\herm,
\end{eqnarray}
where
\begin{eqnarray}
	\bm{A}_{\text{ns},\text{HFR}} = \text{diag}(
		\bm{A}_{h_l,\text{NI}(h_l)}^{\phi, \phi}, 
		\quad
		\underbrace{
		\bm{A}_{h_l-1,1}^{\phi, \phi}, \quad \cdots, \quad \bm{A}_{h_l-1,\Gamma_{h_l-1}}^{\phi, \phi}
		}_{(h_l-1)\text{-th level}},
		\quad
		\cdots,
		\quad
		\underbrace{
		\bm{A}_{h,1}^{\phi, \phi}, \quad \cdots, \quad \bm{A}_{h,\Gamma_h}^{\phi, \phi}
		}_{h\text{-th level}}
		)
\end{eqnarray}
and
\begin{eqnarray}
	\bm{Q}_{\ast, \text{HFR}} = [
	\bm{I},
	\quad
	\underbrace{
	\bm{Q}_{\ast,1,h_l-1,1}, \quad \cdots, \quad \bm{Q}_{\ast,1,h_l-1,\Gamma_{h_l-1}}
	}_{(h_l-1)\text{-th level}},
	\quad
	\cdots,
	\quad
	\underbrace{
	\prod_{\lambda=h_l-1}^h \bm{Q}_{w, 1, \lambda, 1'}, \quad \cdots, \quad
	\prod_{\lambda=h_l-1}^h \bm{Q}_{w, 1, \lambda, \Gamma_h'}
	}_{h\text{-th level}}
	]
\end{eqnarray}
are the sparsified system matrix and transform matrix relative to the 
high frequency regime, respectively.

Substitute Eq. \eqref{eq_ahlphiphi} into the non-standard form system matrix
with respect to the wavelet basis in the low frequency regime, 
we would get the linear system
$\bm{b} = \overline{\bm{Q}_w} \tilde{\bm{A}}_\text{ns} \bm{Q}_\chi^\herm \bm{x}$
for high frequency problems, in which
\begin{eqnarray}
	\begin{split}
		\bm{Q}_\ast = & [ 
		\bm{Q}_{\ast,0,L},\quad
		\bm{Q}_{\ast,1,L},\quad
		\cdots, \quad
		\prod_{l=L}^{h_l+1} \bm{Q}_{\ast,1,l} \cdot \bm{Q}_{\ast,0,h_l}, \quad
		\prod_{l=L}^{h_l+1} \bm{Q}_{\ast,1,l} \cdot \bm{Q}_{\ast, \text{HFR}} ] \\
		= & [
		\bm{Q}_{\ast,0,L},\quad
		\bm{Q}_{\ast,1,L},\quad
		\cdots, \quad
		\prod_{l=L}^{h_l+1} \bm{Q}_{\ast,1,l} \cdot \bm{Q}_{\ast,0,h_l}, \quad
		\prod_{l=L}^{h_l+1} \bm{Q}_{\ast,1,l} \cdot \bm{Q}_{\ast,1,h_l}, \quad \\
		& 
		\underbrace{
		\prod_{l=L}^{h_l} \bm{Q}_{\ast,1,l} \cdot \bm{Q}_{\ast,1,h-1,1}, \quad
		\cdots, \quad
		\prod_{l=L}^{h_l} \bm{Q}_{\ast,1,l} \cdot \prod_{\lambda=h_l-1}^h \bm{Q}_{\ast,1,\lambda,\Gamma_h}
		}_{\text{for all directional cones in the high frequency regime}}
		],
	\end{split}
\end{eqnarray}
and
\begin{eqnarray} \label{eq_ans_curvelet}
	\tilde{\bm{A}}_\text{ns} = \left[ \begin{matrix}
		\bm{A}_L^{\psi, \psi} & \bm{A}_L^{\psi, \phi} & ~ & ~ & ~ & ~ & ~ & ~ & ~ & ~ \\
		\bm{A}_L^{\phi, \psi} & ~ & ~ & ~ & ~ & ~ & ~ & ~ & ~ & ~ \\
		~ & ~ & \ddots & ~ & ~ & ~ & ~ & ~ & ~ & ~ \\
		~ & ~ & ~ & \bm{A}_{h_l+1}^{\psi, \psi} & \bm{A}_{h_l+1}^{\psi, \phi} & ~ & ~ & ~ & ~ & ~ \\
		~ & ~ & ~ & \bm{A}_{h_l+1}^{\phi, \psi} & ~ & ~ & ~ & ~ & ~ & ~ \\
		~ & ~ & ~ & ~ & ~ & \bm{A}_{h_l}^{\psi, \psi} & \bm{A}_{h_l}^{\psi, \phi} & ~ & ~ & ~ \\
		~ & ~ & ~ & ~ & ~ & \bm{A}_{h_l}^{\phi, \psi} & \bm{A}_{h_l, \text{NI}(h_l)}^{\phi, \phi} & ~ & ~ & ~ \\
		~ & ~ & ~ & ~ & ~ & ~ & ~ & \bm{A}_{h_l-1,1}^{\phi, \phi} & ~ & ~ \\
		~ & ~ & ~ & ~ & ~ & ~ & ~ & ~ & \ddots & ~ \\
		~ & ~ & ~ & ~ & ~ & ~ & ~ & ~ & ~ & \bm{A}_{h,\Gamma_h}^{\phi, \phi}
	\end{matrix} \right]
\end{eqnarray}
is the explicitly-sparse representation of the system matrix in the new 
non-standard form for high frequency cases.

The matrix blocks in $\tilde{\bm{A}}_\text{ns}$ for the low frequency regime
is computed in the same manner as that in WBM, and $\bm{A}_{l, \gamma}^{\phi, \phi}$ 
for the high frequency regime can be
computed efficiently with the M2L translation.
It is worthnoting that each submatrix block in $\bm{A}_{l, \gamma}^{\phi, \phi}$ 
corresponds to an M2L operation, 
thus there are only $O(N \log N)$ blocks. 
The size of each block is of order $O(1)$ due to the directional low rank 
property of the kernel. Therefore, there are only $O(N \log N)$ nonzero elements
in $\tilde{\bm{A}}_\text{ns}$, and the computational complexity
of the matrix computation is also $O(N \log N)$.
Then the computational cost of $\bm{b = Ax}$ can be brought down to log-linear 
by computing fast curvelet transform 
$\tilde{\bm{x}} = \bm{Q}_\chi^\herm \bm{x}$, 
matrix vector multiplication in the curvelet space
$\tilde{\bm{b}} = \tilde{\bm{A}}_\text{ns} \tilde{\bm{x}}$, 
and the inverse fast curvelet transform
$\bm{b} = \overline{\bm{Q}_w} \tilde{\bm{b}}$ in sequence.

\section{A-posteriori compression}
\label{sec_aposteriori}

In the study of WBM, it is found that
that although the system matrix is sparsified with wavelets for low
frequency cases, there are still many tiny elements. The memory cost can be 
significantly reduced by leaving out these tiny elements, which is called
the \emph{a-posteriori} compression technique
\cite{xiaojinyou2009aposteriori}. Considering the relationship between our CBM 
and the WBM, it is reasonable to expect that there would be also many discardable
elements in the explicitly-sparse representation obtained in Section 
\ref{sec_curvelet}, and its sparsify can be enhanced by the \emph{a-posteriori}
compression technique.

It is noteworthy that our aim is to develop a nearly optimize explicitly-sparse
representation of the system matrix with controllable accuracy.
However, the primitive \emph{a-posteriori} compression technique is developed for
fast boundary element analysis and the threshold is designed to preserve the
convergence rate of the boundary element method with piecewise constant elements.
Thus the threshold for leaving out should be re-defined in our CBM.

For each matrix block $\bm{B}$ in $\tilde{\bm{A}}_\text{ns}$, it is compressed by
\begin{eqnarray}
	\tilde{B}_{ij} = \begin{cases}
		B_{ij}, & \quad \| B_{ij} \| \ge \eta, \\
		0, & \quad \| B_{ij} \| < \eta.
	\end{cases}
\end{eqnarray}
The induced relative error would be controllable when
\begin{eqnarray}
	\frac{\| \bm{B} - \tilde{\bm{B}} \|_1}{\| \bm{B} \|_1} 
	\le \eta m / \| \bm{B} \|_1 = \varepsilon,
\end{eqnarray}
where $L_1$ norm is used in the error estimation since it is easy to evaluate.
Therefore, the threshold may be taken as
\begin{eqnarray} \label{eq_threshold_aposteriori}
	\eta = \varepsilon  \| \bm{B} \|_1 / m.
\end{eqnarray}

The compressed matrix block $\tilde{\bm{B}}$ should be stored in sparse mode
to save memory. For each nonzero element, it is stored with a complex variable for
its value and an interger for its position. Thus $\tilde{\bm{B}}$ with $N_\text{nz}$
nonzeros can be stored with $N_\text{nz} [ \texttt{sizeof(int)} + 
\texttt{sizeof(complex)} ]$ memory. Notice that the matrix block before a-posteriori 
compression requires $mn \cdot \texttt{sizeof(complex)}$ memory, thus the 
a-posteriori compression should be only carried out for $\bm{B}$ when there are 
\begin{eqnarray}
	N_\text{nz} < mn \frac{\texttt{sizeof(complex)}}{\texttt{sizeof(int)} + 
	\texttt{sizeof(complex)}}
\end{eqnarray}
undiscardable elements.

\section{Numerical results}
\label{sec_numres}

The CBM implemented in this work is transformed from the optimized dFMM
\cite{caoyanchuang2015fdbem}, that is, the low rank approximation are constructed
with equivalent densities, and compressed S2M, M2M, M2L, L2L and L2T translation
matrices in \cite{caoyanchuang2015fdbem} are used in the curvelet construction
and matrix computation.

%
%
%

The program is implemented serially in C++.
Its performance is demonstrated with a group of numerical examples
in this section. All the examples are carried out with a Chinese
FT2000+ CPU (2.2 GHz) and 128 GB RAM.

\subsection{Accuracy and computational complexity}

First we study the accuracy and computational complexity of our CBM with 
the summation \eqref{eq_nbody} and single layer Helmholtz kernel 
\eqref{eq_greenfunc}.
The points $\{ \bm{x}_i \}$ and 
$\{ \bm{y}_j \}$ are sampled on the surface of a unit sphere with about
10 points per wavelength, and the densities $\{ \sigma_j \}$ are randomly 
defined. The relative error $\varepsilon_\text{a}$ of our CBM is estimated with
\begin{eqnarray}
	\varepsilon_\text{a} = \sqrt{
		\frac{ \sum_{i=1}^{N_\text{t}}
		\left\| f_i^\text{(a)} - f_i^\text{(d)}\right\|_2^2 }
		{ \sum_{i=1}^{N_\text{t}} \left\| f_i^\text{(d)} \right\|_2^2 }
		},
\end{eqnarray}
where $\{ f_i^\text{(a)} \}$ are the potentials on $N_\text{t}$ randomly 
selected points computed with our CBM, and $\{ f_i^\text{(d)} \}$ are the 
potentials evaluated by direct summation. In this work, $N_\text{t}=500$
points are selected in the error estimation.

\begin{figure}[ht]
	\centering
	\includegraphics[width=0.4\textwidth]{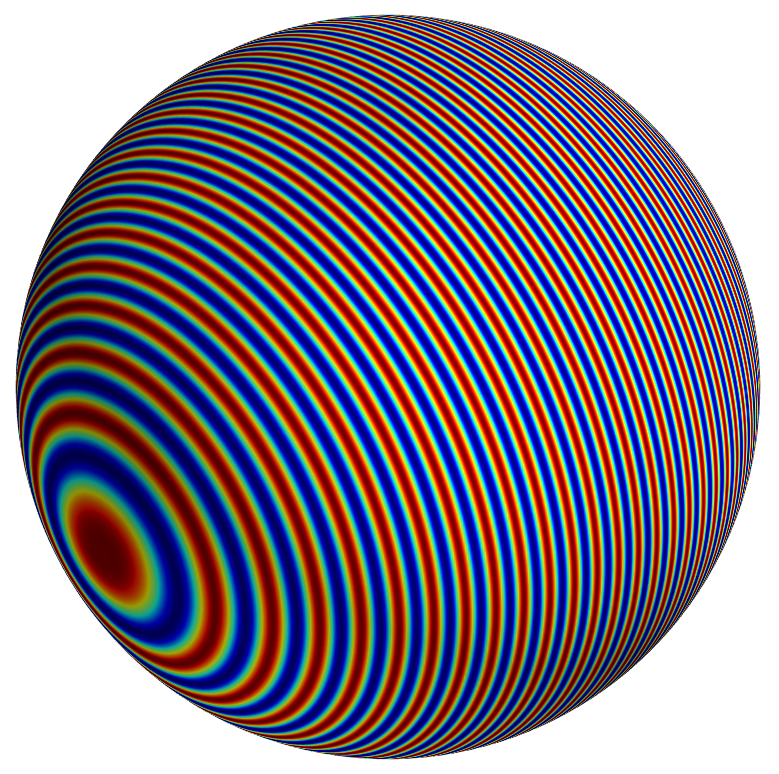}
	\caption{A sphere with $\kappa D = 201.1$.}
	\label{fig_sphere}
\end{figure}

The summation is computed multiple times with various choices of 
$\varepsilon$ and dimensionless diameters $\kappa D$.
In the highest frequency case with $\kappa D = 201.1$, the diameter is 
32 wavelengths, as illustrated in Figure \ref{fig_sphere},
and over 1 million points are sampled on the boundary.
The computational results are summarized in Table \ref{tab_single},
where $T_\text{c}, T_\text{m}, T_\text{p}, T_\text{t}$ are the computational time
cost by curvelet construction, sparse system matrix computation, matrix-vector multiplication,
and the overall running time, $M_\text{Q}, M_\text{m}, M_\text{t}$ are the memory
cost by the $\bm{Q}$ matrices, sparse system matrix, and the overall memory consumption,
respectively.

\begin{table} [ht]
	\centering
	\caption{Results of the sphere with the single-layer kernel.}
	\vspace{1em}
	\begin{tabular}{rrr|rrrr|rrr|r}
		\hline
		\multirow{2}{*}{$\varepsilon_0$}	& \multirow{2}{*}{$\kappa D$}	& \multirow{2}{*}{$N$}	& \multicolumn{4}{c|}{Time cost (sec)}	& \multicolumn{3}{c|}{Memory cost (MB)}	& \multirow{2}{*}{$\varepsilon_\text{a}$} \\
		\cline{4-10}
		& & & $T_\text{c}$	& $T_\text{m}$	& $T_\text{p}$	& $T_\text{t}$	& $M_\text{Q}$		& $M_\text{m}$		& $M_\text{t}$		\\
		\hline
		1e-3			&    12.6		&     4,608	&     1					&     5					&$<$1					&      6				&     12				&     51				&      68				& 1.0e-3	\\
		1e-3			&    25.1		&    18,432	&    15					&    27					&$<$1					&     43				&    125				&    294				&     518				& 2.7e-4	\\
		1e-3			&    50.3		&    73,728	&   129					&   103					&   2					&    236				&  1,013				&  1,406				&   2,671				& 4.4e-3	\\
		1e-3			&   100.5		&   294,912	&   812					&   429					&  11					&  1,265				&  6,414				&  6,083				&  14,011				& 5.5e-3	\\
		1e-3			&   201.1		& 1,143,072	& 4,054					& 1,866					&  79					&  6,053				& 31,987				& 26,241				&  62,579				& 8.8e-3	\\
		\hline
		1e-6			&    12.6		&     4,608	&     5					&    44					&$<$1					&     51				&     12				&    208				&     249				& 1.2e-6	\\
		1e-6			&    25.1		&    18,432	&   173					&   351					&   1					&    528				&    518				&  2,641				&   3,465				& 6.2e-6	\\
		1e-6			&    50.3		&    73,728	& 1,796					& 1,200					&   8					&  3,009				&  5,465				&  9,826				&  16,923				& 8.3e-6	\\
		1e-6			&   100.5		&   294,912	&14,238					& 5,074					&  83					& 19,413				& 38,089				& 46,476				&  94,412				& 9.5e-6	\\
		\hline
	\end{tabular}
	\label{tab_single}
\end{table}

Generally the overall error $\varepsilon_\text{a}$ is of the same magnitude with 
the predetermined parameter $\varepsilon$.
The error $\varepsilon_\text{a}$ seems to grow with 
$\log N$. This is reasonable since in each level, some error
of order $O(\varepsilon)$ is introduced, and there are $O(\log N)$ levels 
in the octree.
Compared with the error of the dFMM in \cite{caoyanchuang2015fdbem} from which
our algorithm is transformed, the error becomes greater. This is because in our
CBM, besides the low rank approximation of the kernel, extra error is 
introduced in the construction of curvelets with quasi-vanishing moments and
the \emph{a-posteriori} compression. 
Nevertheless, the overall error is still controllable.

In the $N$-body problem evaluations, most of the computational time are cost
by construction of the explicitly-sparse representation of the system matrix, 
including the curvelet construction and matrix computation.
Once the sparse representation is obtained, potentials on the target points
can be evaluated very efficiently. 
A majority of the memory cost are taken by the storage of $\bm{Q}$ matrices and
the sparse system matrix $\tilde{\bm{A}}_\text{ns}$, which are necessary
for the curvelet transformation and matrix-vector multiplication
in the curvelet space.

\begin{figure}[ht]
	\centering
	\subfigure[Overall time cost.]{
		\includegraphics[width=0.47\textwidth]{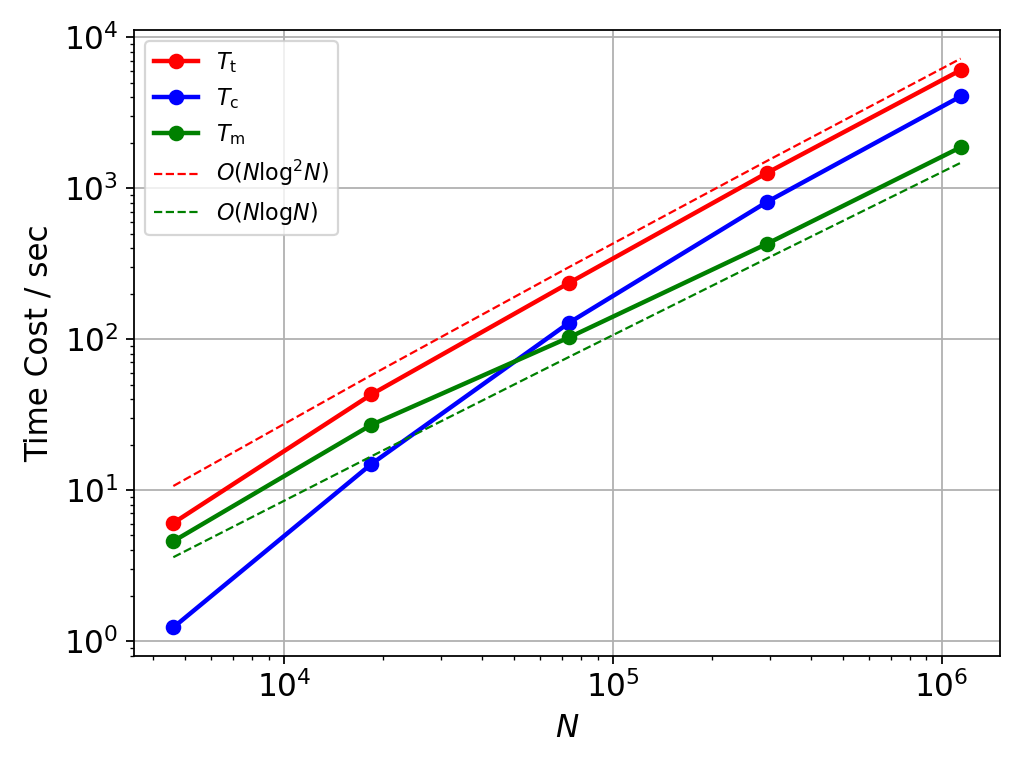}
		\label{fig_single_time}
	}
	\subfigure[Overall memory consumption.]{
		\includegraphics[width=0.47\textwidth]{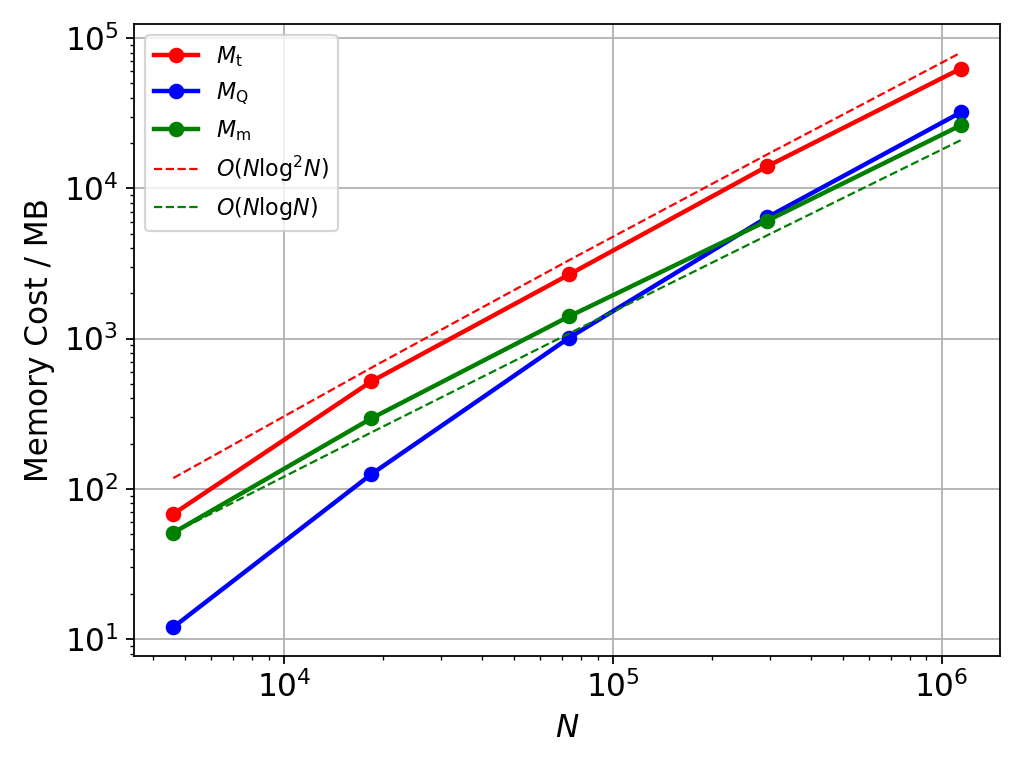}
	}
	\caption{Computational cost of the summation on a unit sphere with the single layer kernel.}
	\label{fig_single}
\end{figure}

The time and memory cost by evaluations with $\varepsilon_0 =$ 1e-3 are 
plotted in Figure \ref{fig_single} to show the computational complexity. 
It shows that the memory taken by storing $\tilde{\bm{A}}_\text{ns}$
is of order $O(N \log N)$, which means there are only log-linear nonzeros
in the explicitly-sparse representation of the system matrix obtained by
our algorithm. The time cost taken by the matrix computation also increases 
at the speed of $O(N \log N)$.

Theoretically, the computational complexity of curvelet construction 
and the overall computation should also be $O(N \log N)$.
For each level with the cube width being $w$, there are
at most $O(D^2/w^2)$ non-empty cubes since the points are sampled on a surface.
And there are at most $O(\kappa^2 w^2)$ directional cones when the cube lies
in the high frequency regime. The computational cost of curvelet construction 
for each directional cone is of order $O(1)$ since the dimension of the moment
matrices and M2M matrices are independent of $\kappa w$. 
Therefore, the computational cost of curvelet construction in each level
is of order $O(D^2 / w^2 \cdot \kappa^2 w^2) = O(\kappa^2 D^2) = O(N)$, 
and the overall computational complexity of curvelet construction in all levels
should be $O(N \log N)$ since there are $O(\log N)$ levels in the octree.

However, it seems the time and memory cost in curvelet construction
increases faster than the theoretical prediction, and 
the overall computational complexity seems $O(N \log^2 N)$.
To figure out the reason, we count the actual number of nonempty 
directional cones in the numerical cases.
Notice that for a cube in the high frequency, if its interaction field
in the directional cone is empty but there are nonempty cubes in the far field
in this cone, the cone is still considered as nonempty, and curvelets and 
scaling functions have to be constructed since they are later required in the 
construction of its parent in higher levels. The number of nonempty cones
differs from cube to cube, thus we compute its average value for each level.
For each level with cube width in term of wavelength $w/\lambda = 
\kappa w / 2\pi$, the average number of nonempty cones $\tilde{\Gamma}_l$ 
for each cone and its theoretical upper bound $\bar{\Gamma}_l$ are listed in 
Table \ref{tab_num_cones}.

\begin{table} [ht]
	\centering
	\caption{Number of nonempty directional cones in the $N$-body problem with
	points on a spherical surface.}
	\vspace{1em}
	\begin{tabular}{r|rrrr|r}
		\hline
		\multirow{2}{*}{$w/\lambda$}	& \multicolumn{4}{c|}{$\tilde{\Gamma}_l$}	& \multirow{2}{*}{$\bar{\Gamma}_l$} \\
		\cline{2-5}
		& $\kappa D = 8\pi$	& $\kappa D = 16 \pi$	& $\kappa D = 32 \pi$	& $\kappa D = 64\pi$	\\
		\hline
		1		& 12.9	& 35.4	& 49.1	&  55.7		&  96	\\
		2		& -		& -		& 69.1	& 142.0		& 384	\\
		\hline
	\end{tabular}
	\label{tab_num_cones}
\end{table}

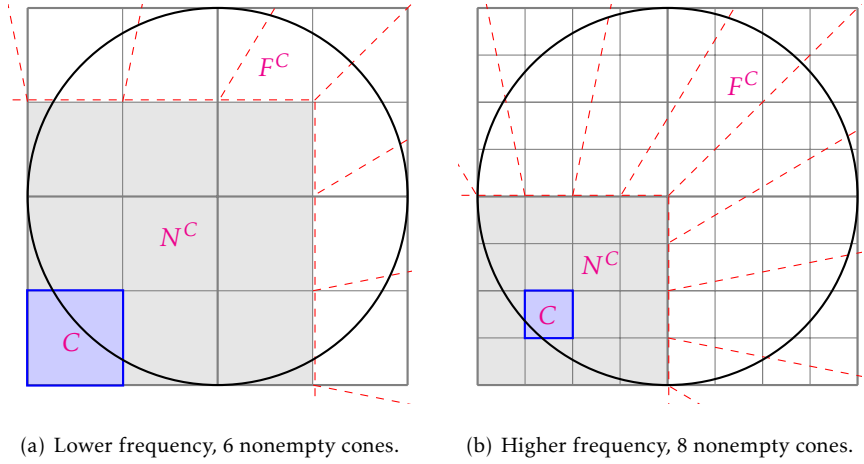
\begin{figure}[ht]
	\centering
	\subfigure[Lower frequency, 6 nonempty cones.] {
		\begin{tikzpicture}[>=latex, scale=2.5]
			\draw[gray!20, thin, fill=gray!20] (0,0) rectangle (1.5,1.5);
			\draw[magenta] (0.80,0.80) node {$N^C$};
			\draw[blue!20, thin, fill=blue!20] (0,0) rectangle (0.5,0.5);
			\draw[magenta] (0.23,0.23) node {$C$};
			\draw[magenta] (1.30,1.70) node {$F^C$};
			\draw[step=1cm, gray, thick] (0,0) grid(2,2);
			\draw[step=0.5cm, gray] (0,0) grid (2,2);
			\clip (-0.1,-0.1) rectangle (2.02,2.02);
			\begin{scope}[xshift=0.25cm, yshift=0.25cm, scale=0.505]
				\draw[blue, thick] (-0.5,-0.5) rectangle (0.5,0.5);
				\draw[red, dashed] (-2.5,-2.5) rectangle (2.5,2.5);
				\draw[red, dashed] (101.3:2.53) -- (101.3:5);
				\draw[red, dashed] (78.7:2.53) -- (78.7:5);
				\draw[red, dashed] (59.0:2.95) -- (59.0:5);
				\draw[red, dashed] (45.0:3.53) -- (45.0:5);
				\draw[red, dashed] (31.0:2.95) -- (31.0:5);
				\draw[red, dashed] (11.3:2.53) -- (11.3:5);
				\draw[red, dashed] (-11.3:2.53) -- (-11.3:5);
			\end{scope}
			\draw[thick] (1,1) circle (1cm);
		\end{tikzpicture}
	}
	\quad
	\subfigure[Higher frequency, 8 nonempty cones.] {
		\begin{tikzpicture}[>=latex, scale=2.5]
			\draw[gray!20, thin, fill=gray!20] (0,0) rectangle (1.0,1.0);
			\draw[magenta] (0.65,0.65) node {$N^C$};
			\draw[blue!20, thin, fill=blue!20] (0.25,0.25) rectangle (0.5,0.5);
			\draw[magenta] (0.37,0.37) node {$C$};
			\draw[magenta] (1.40,1.60) node {$F^C$};
			\draw[step=1cm, gray, thick] (0,0) grid(2,2);
			\draw[step=0.5cm, gray] (0,0) grid (2,2);
			\draw[step=0.25cm, gray, thin] (0,0) grid (2,2);
			\clip (-0.1,-0.1) rectangle (2.02,2.02);
			\begin{scope}[xshift=0.375cm, yshift=0.375cm, scale=0.2525]
				\draw[blue, thick] (-0.5,-0.5) rectangle (0.5,0.5);
				\draw[red, dashed] (-2.5,-2.5) rectangle (2.5,2.5);
				\draw[red, dashed] (121.0:2.95) -- (121.0:10);
				\draw[red, dashed] (101.3:2.53) -- (101.3:10);
				\draw[red, dashed] (78.7:2.53) -- (78.7:10);
				\draw[red, dashed] (59.0:2.95) -- (59.0:10);
				\draw[red, dashed] (45.0:3.53) -- (45.0:10);
				\draw[red, dashed] (31.0:2.95) -- (31.0:10);
				\draw[red, dashed] (11.3:2.53) -- (11.3:10);
				\draw[red, dashed] (-11.3:2.53) -- (-11.3:10);
				\draw[red, dashed] (-31.0:2.95) -- (-31.0:10);
			\end{scope}
			\draw[thick] (1,1) circle (1cm);
		\end{tikzpicture}
	}
	\caption{Increasement of the number of nonempty directional cones
	with frequency for the same $w/\lambda$.
	Illustrated with a two-dimensional case.
	Three-dimensional cases would be similar.
	}
	\label{fig_curvature}
\end{figure}

It is shown that in real-world numerical cases, the number of nonempty
directional cones may be much less than the upper bound $\bar{\Gamma}_l \sim 
O(\kappa^2 w^2)$, and it tends to increase with frequency.
This may be because the dimensionless curvature $\lambda/R$
of the surface on which the points locate reduces with the frequency. 
When the frequency is not so high, the curvature is relatively great,
and the far-field surface tends to locate in a small number directional cones.
But when the frequency gets higher, the curvature decreases and the surface
tends to occupy more directional cones, as illustrated in Figure \ref{fig_curvature}.
Nevertheless, the number of directional cones will never exceed $\bar{\Gamma}_l
\sim O(\kappa^2 w^2)$. Thus the computational complexity of curvelet construction
should approach $O(N \log N)$ asymptotically as the frequency keep increasing.
This trend has been partially demonstrated in Figure \ref{fig_single},
in which the slope of $T_\text{c}$ reduces as the frequency increases.
It is expected that the slope would further reduce to approximately 
parallel with that of $O(N \log N)$ if results of higher frequency cases are 
provided. However, they are not computed in this work due to the memory limit 
of our computer.

\subsection{Kernels of other layers}

The performance of our CBM with double, adjoint and quadrapole layers are 
studied in this section. The points are also sampled on the surface of a unit 
sphere, and their normals points outwards.
The numerical results are listed in Table
\ref{tab_double}--\ref{tab_quadrapole}.

\begin{table} [ht]
	\centering
	\caption{Results of the sphere with the double-layer kernel.}
	\vspace{1em}
	\begin{tabular}{rrr|rrrr|rrr|r}
		\hline
		\multirow{2}{*}{$\varepsilon_0$}	& \multirow{2}{*}{$\kappa D$}	& \multirow{2}{*}{$N$}	& \multicolumn{4}{c|}{Time cost (sec)}	& \multicolumn{3}{c|}{Memory cost (MB)}	& \multirow{2}{*}{$\varepsilon_\text{a}$} \\
		\cline{4-10}
		& & & $T_\text{c}$	& $T_\text{m}$	& $T_\text{p}$	& $T_\text{t}$	& $M_\text{Q}$		& $M_\text{m}$		& $M_\text{t}$		\\
		\hline
		1e-3			&    12.6		&     4,608	&     1					&     5					&$<$1					&      6				&     12				&     61				&      77				& 5.0e-4	\\
		1e-3			&    25.1		&    18,432	&    15					&    26					&$<$1					&     43				&    125				&    388				&     559				& 1.1e-3	\\
		1e-3			&    50.3		&    73,728	&   124					&    95					&   2					&    225				&  1,013				&  1,901				&   3,177				& 1.2e-3	\\
		1e-3			&   100.5		&   294,912	&   786					&   403					&  35					&  1,238				&  6,465				&  8,515				&  16,518				& 1.7e-3	\\
		1e-3			&   201.1		& 1,143,072	& 4,065					& 1,738					&  83					&  5,941				& 32,884				& 38,000				&  74,970				& 1.8e-3	\\
		\hline                           
		1e-6			&    12.6		&     4,608	&     6					&    48					&$<$1					&     55				&      6				&    316				&     351				& 5.8e-7	\\
		1e-6			&    25.1		&    18,432	&   164					&   358					&   1					&    532				&    531				&  2,971				&   3,807				& 1.4e-6	\\
		1e-6			&    50.3		&    73,728	& 1,869					& 1,241					&   9					&  3,124				&  5,854				& 12,161				&  19,753				& 1.5e-6	\\
		1e-6			&   100.5		&   294,912	&15,357					& 5,252					&  80					& 20,708				& 40,368				& 56,058				& 106,432				& 1.7e-6	\\
		\hline
	\end{tabular}
	\label{tab_double}
\end{table}

\begin{table} [ht]
	\centering
	\caption{Results of the sphere with the adjoint-layer kernel.}
	\vspace{1em}
	\begin{tabular}{rrr|rrrr|rrr|r}
		\hline
		\multirow{2}{*}{$\varepsilon_0$}	& \multirow{2}{*}{$\kappa D$}	& \multirow{2}{*}{$N$}	& \multicolumn{4}{c|}{Time cost (sec)}	& \multicolumn{3}{c|}{Memory cost (MB)}	& \multirow{2}{*}{$\varepsilon_\text{a}$} \\
		\cline{4-10}
		& & & $T_\text{c}$	& $T_\text{m}$	& $T_\text{p}$	& $T_\text{t}$	& $M_\text{Q}$		& $M_\text{m}$		& $M_\text{t}$		\\
		\hline
		1e-3			&    12.6		&     4,608		&     1				&     5					&$<$1					&      6				&     12				&     60				&      77				& 5.0e-4	\\
		1e-3			&    25.1		&    18,432		&    15				&    26					&$<$1					&     41				&    125				&    388				&     559				& 1.1e-3	\\
		1e-3			&    50.3		&    73,728		&   127				&    96					&   2					&    228				&  1,012				&  1,901				&   3,177				& 1.2e-3	\\
		1e-3			&   100.5		&   294,912		&   786				&   403					&  34					&  1,237				&  6,465				&  8,515				&  16,518				& 1.7e-3	\\
		1e-3			&   201.1		& 1,143,072		& 4,068				& 1,732					&  82					&  5,937				& 32,884				& 38,000				&  74,970				& 1.8e-3	\\
		\hline                          
		1e-6			&    12.6		&     4,608		&     5				&    44					&$<$1					&     51				&      6				&    316				&     351				& 5.8e-7	\\
		1e-6			&    25.1		&    18,432		&   160				&   352					&   1					&    515				&    531				&  2,971				&   3,807				& 1.4e-6	\\
		1e-6			&    50.3		&    73,728		& 1,875				& 1,234					&  10					&  3,124				&  5,854				& 12,161				&  19,753				& 1.5e-6	\\
		1e-6			&   100.5		&   294,912		&15,494				& 5,254					&  80					& 20,847				& 40,368				& 56,058				& 106,432				& 1.7e-6	\\
		\hline
	\end{tabular}
	\label{tab_adjoint}
\end{table}

\begin{table} [ht]
	\centering
	\caption{Results of the sphere with the quadrapole-layer kernel.}
	\vspace{1em}
	\begin{tabular}{rrr|rrrr|rrr|r}
		\hline
		\multirow{2}{*}{$\varepsilon_0$}	& \multirow{2}{*}{$\kappa D$}	& \multirow{2}{*}{$N$}	& \multicolumn{4}{c|}{Time cost (sec)}	& \multicolumn{3}{c|}{Memory cost (MB)}	& \multirow{2}{*}{$\varepsilon_\text{a}$} \\
		\cline{4-10}
		& & & $T_\text{c}$	& $T_\text{m}$	& $T_\text{p}$	& $T_\text{t}$	& $M_\text{Q}$		& $M_\text{m}$		& $M_\text{t}$		\\
		\hline
		1e-3			&    12.6		&     4,608	&     1					&     5					&$<$1					&      6				&     12				&     46				&      62				& 8.0e-4	\\
		1e-3			&    25.1		&    18,432	&    14					&    25					&$<$1					&     41				&    125				&    159				&     456				& 2.2e-4	\\
		1e-3			&    50.3		&    73,728	&   120					&    94					&   2					&    219				&  1,023				&    722				&   2,382				& 3.7e-3	\\
		1e-3			&   100.5		&   294,912	&   773					&   394					&  10					&  1,192				&  6,380				&  3,562				&  12,533				& 4.5e-3	\\
		1e-3			&   201.1		& 1,143,072	& 3,961					& 1,712					&  79					&  5,808				& 33,159				& 17,018				&  58,595				& 5.7e-3	\\
		\hline                          
		1e-6			&    12.6		&     4,608	&     6					&    47					&$<$1					&     54				&      0				&    324				&     353				& 4.0e-7	\\
		1e-6			&    25.1		&    18,432	&   164					&   352					&   1					&    519				&    546				&  2,347				&   3,196				& 6.4e-6	\\
		1e-6			&    50.3		&    73,728	& 1,943					& 1,241					&   9					&  3,196				&  6,216				&  8,624				&  16,677				& 6.5e-6	\\
		1e-6			&   100.5		&   294,912	&17,068					& 5,340					&  77					& 22,502				& 43,356				& 39,995				&  94,160				& 1.2e-5	\\
		\hline
	\end{tabular}
	\label{tab_quadrapole}
\end{table}

It is shown that the performance for the adjoint layer kernel is approximately
the same with that for the double layer. This is reasonable since their matrices
are in fact the transpose of each other, i.e., $\bm{A}^\text{(a)} = 
\overline{\bm{A}^\text{(d)}}$, where the superscript (a) and (d) represents the 
adjoint and double layer kernel, respectively. Since in our CBM we get the sparse 
representation $\bm{A} = \overline{\bm{Q}_w} \tilde{\bm{A}}_\text{ns} \bm{Q}_\chi$,
thus theoretically, $\bm{Q}_w^\text{(a)} = \bm{Q}_\chi^\text{(d)}, 
\tilde{\bm{A}}_\text{ns}^\text{(a)} = 
\overline{\tilde{\bm{A}}_\text{ns}^\text{(d)}}$, and 
$\bm{Q}_\chi^\text{(a)} = \bm{Q}_w^\text{(d)}$. 
Hence the computations for the adjoint layer kernel are the same with that
for the double layer kernel, except that the position of target and source
points are exchanged in the algorithm.
The slight differences in the computational cost may come from the random 
choice of equivalent points in construction of low rank approximation for 
the kernel in high frequency regimes \cite{engquist2007fda, engquist2010dfmm}.

The memory cost $M_\text{m}$ by storing $\tilde{A}_\text{ns}$ for the 
quadrapole layer is much less than that for the single layer.
This is because, for the quadrapole layer kernel
\begin{eqnarray}
	\frac{\partial^2 G(\bm{x}, \bm{y})}{\partial \bm{n_x} \partial \bm{n_y}}
   	= \left[(k^2 r^2 + 3ikr - 3) (\bm{\hat{r}} \cdot \bm{\hat{n}_x}) 
	(\bm{\hat{r}} \cdot \bm{\hat{n}_y}) + (1-ikr)(\bm{\hat{n}_x} \cdot
	\bm{\hat{n}_y})\right] \frac{\text{e}^{ikr}}{4\pi r^3},
\end{eqnarray}
its hypersingularity makes it increases sharply
when $\bm{x}$ approaches $\bm{y}$. Thus the elements in the near field
can be much greater than that in the case with the single layer kernel.
The norm of the transformed matrix block $\bm{A}^{\phi,\phi}$,
and the threshold for the \emph{a-posteriori} compression 
\eqref{eq_threshold_aposteriori} would also become larger. Consequently,
more tiny elements could be discarded in the \emph{a-posteriori} compression,
leaving less nonzero elements in the final explicitly-sparse representation,
and less memory is consumed.
The overall error $\varepsilon_\text{a}$ becomes greater 
but is still approximately of the same order with $\varepsilon$.

The memory cost $M_\text{m}$ for the double and adjoint layer, however, is
greater than that for the single layer.
This is because, for the double layer kernel
\begin{eqnarray}
	\frac{\partial G(\bm{x}, \bm{y})}{\partial \bm{n_y}} = (1 - ikr)
		\cdot (\bm{\hat{r}} \cdot \bm{\hat{n}_y})\frac{\text{e}^{ikr}}{4\pi r^2}
\end{eqnarray}
and the adjoint kernel
\begin{eqnarray}
	\frac{\partial G(\bm{x}, \bm{y})}{\partial \bm{n_x}} = -(1 - ikr)
		\cdot (\bm{\hat{r}} \cdot \bm{\hat{n}_x})\frac{\text{e}^{ikr}}{4\pi r^2},
\end{eqnarray}
the $\bm{\hat{r}} \cdot \bm{\hat{n}}$ term makes them approaches 0 when $\bm{x}$
approaches $\bm{y}$ since $\bm{x}$ and $\bm{y}$ are sampled on the surface
with $\bm{n}$ is the surface normal. This makes the elements in the near field
becomes much smaller than that in the case with the single layer kernel, and 
results in a smaller threhold for the \emph{a-posteriori} compression.
Consequently, more nonzero elements are left in the explicitly-sparse
representation and more memory is consumed.
The overall error $\varepsilon_\text{a}$ becomes lower probably because
less error is introduced in the \emph{a-posteriori} compression.

\subsection{Different geometries}

\begin{figure}[ht]
	\centering
	\includegraphics[width=0.6\textwidth]{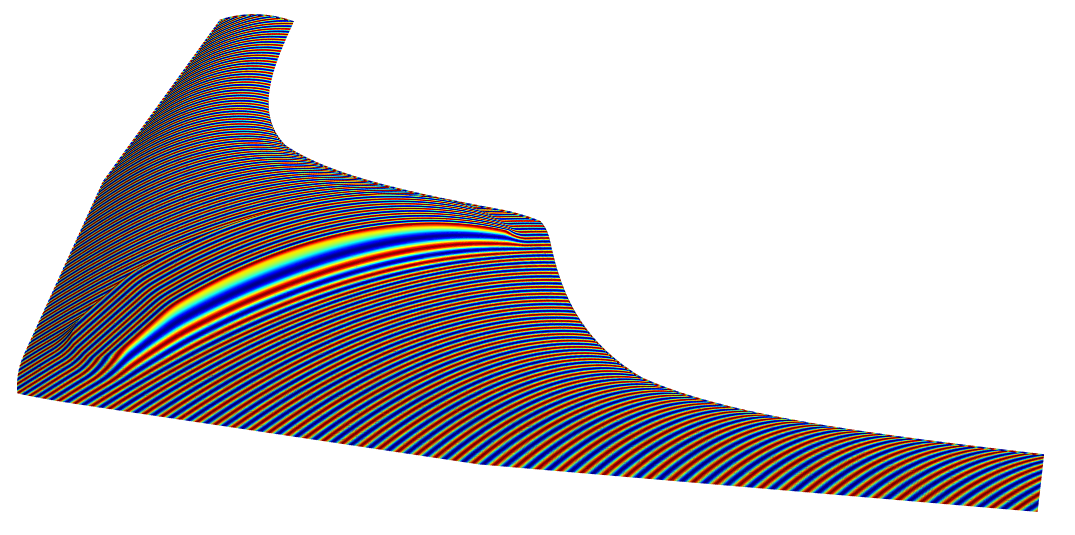}
	\caption{A aircraft with $\kappa D = 603.2$.}
	\label{fig_x45x}
\end{figure}

\begin{figure}[ht]
	\centering
	\includegraphics[width=0.7\textwidth]{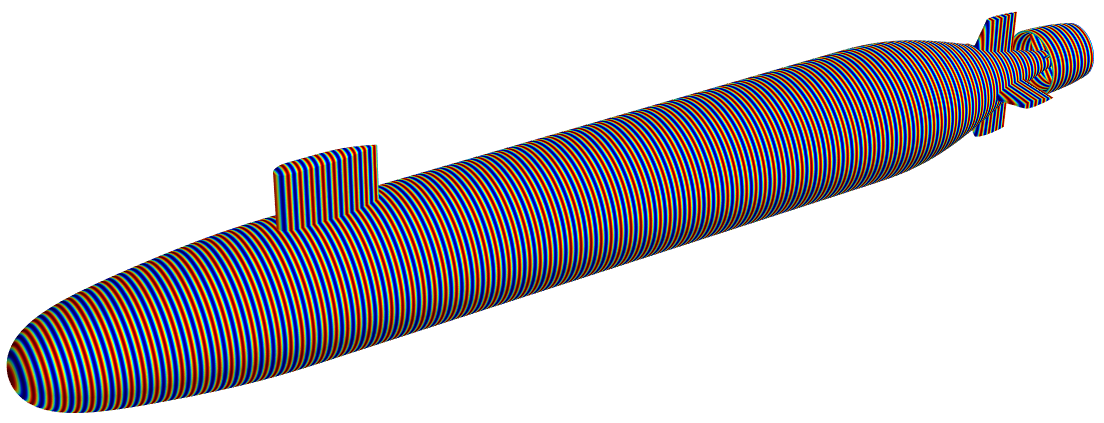}
	\caption{A submarine with $\kappa D = 684.3$.}
	\label{fig_suboff}
\end{figure}

An aircraft X45X and a submarine DARPA suboff is used to study the performance 
of our CBM on different geometries.
The CAD models are downloaded from \emph{grabcab.com}.
Various frequencies and different choices of $\varepsilon$ are considered.
In the highest frequency case, the size of the aircraft is about 96 wavelengths, 
and the submarine is about 109 wavelengths long, as depicted in Figure
\ref{fig_x45x} and \ref{fig_suboff}.
The points are also sampled with approximately 10 points per wavelength. 
The summations are computed with the single layer kernel.
The computational results are listed in Table \ref{tab_x45x} and \ref{tab_suboff},
in which the log-linear complexity is shown.
Compared with summations on a sphere, higher frequency cases are computed 
for the aircraft and submarine within approximately the same computational cost,
showing that our CBM is more efficient for flattened and elongated geometries.
This is because with such geometries, more directional cones are empty, and 
less directional cones has to be addressed.

\begin{table} [ht]
	\centering
	\caption{Results of the aircraft with the single-layer kernel.}
	\vspace{1em}
	\begin{tabular}{rrr|rrrr|rrr|r}
		\hline
		\multirow{2}{*}{$\varepsilon_0$}	& \multirow{2}{*}{$\kappa D$}	& \multirow{2}{*}{$N$}	& \multicolumn{4}{c|}{Time cost (sec)}	& \multicolumn{3}{c|}{Memory cost (MB)}	& \multirow{2}{*}{$\varepsilon_\text{a}$} \\
		\cline{4-10}
		& & & $T_\text{c}$	& $T_\text{m}$	& $T_\text{p}$	& $T_\text{t}$	& $M_\text{Q}$		& $M_\text{m}$		& $M_\text{t}$		\\
		\hline
		1e-3			&     75.4		&    22,992	&    17					&    26					&$<$1					&     45				&    155				&    308				&     565				& 2.0e-3	\\
		1e-3			&    150.8		&    88,608	&   114					&   107					&   2					&    226				&  1,026				&  1,738				&   3,024				& 1.4e-3	\\
		1e-3			&    301.6		&   350,388	&   763					&   501					&  25					&  1,305				&  6,547				&  8,180				&  16,121				& 1.7e-3	\\
		1e-3			&    603.2		& 1,401,348	& 3,920					& 2,035					&  99					&  6,139				& 35,612				& 34,059				&  74,600				& 1.7e-3	\\
		\hline
		1e-6			&     75.4		&    22,992	&   149					&   283					&   2					&    438				&    614				&  2,557				&   3,454				& 9.5e-7	\\
		1e-6			&    150.8		&    88,608	& 1,091					& 1,370					&   9					&  2,481				&  4,708				& 13,770				&  20,157				& 8.9e-7	\\
		1e-6			&    301.6		&   350,388	& 6,570					& 5,910					&  82					& 12,603				& 30,559				& 64,361				& 102,260				& 1.2e-6	\\
		\hline
	\end{tabular}
	\label{tab_x45x}
\end{table}

\begin{table} [ht]
	\centering
	\caption{Results of the submarine with the single-layer kernel.}
	\vspace{1em}
	\begin{tabular}{rrr|rrrr|rrr|r}
		\hline
		\multirow{2}{*}{$\varepsilon_0$}	& \multirow{2}{*}{$\kappa D$}	& \multirow{2}{*}{$N$}	& \multicolumn{4}{c|}{Time cost (sec)}	& \multicolumn{3}{c|}{Memory cost (MB)}	& \multirow{2}{*}{$\varepsilon_\text{a}$} \\
		\cline{4-10}
		& & & $T_\text{c}$	& $T_\text{m}$	& $T_\text{p}$	& $T_\text{t}$	& $M_\text{Q}$		& $M_\text{m}$		& $M_\text{t}$		\\
		\hline
		1e-3			&  85.5 		&    21,348	&    12					&    25					&$<$1					&     38				&     98				&    301				&     498				& 6.8e-4	\\
		1e-3			& 171.1 		&    84,204	&   127					&   111					&   2					&    243				&  1,060				&  1,556				&   2,887				& 1.0e-3	\\
		1e-3			& 342.1 		&   330,600	&   770					&   459					&  13					&  1,258				&  6,493				&  7,349				&  15,477				& 1.1e-3	\\
		1e-3			& 684.3 		& 1,342,260	& 4,435					& 1,911					&  98					&  6,527				& 39,701				& 29,105				&  76,253				& 2.4e-3	\\
		\hline
		1e-6			&  85.5 		&    21,348	&    92					&   238					&   1					&    335				&    327				&  2,179				&   2,668				& 7.4e-7	\\
		1e-6			& 171.1 		&    84,204	& 1,101					& 1,288					&   9					&  2,407				&  4,818				& 12,903				&  19,489				& 1.0e-6	\\
		1e-6			& 342.1 		&   330,600	& 7,018					& 5,354					&  79					& 12,503				& 32,125				& 59,739				& 101,496				& 1.1e-6	\\
		\hline
	\end{tabular}
	\label{tab_suboff}
\end{table}

\section{Conclusion}

A nearly optimal explicitly-sparse representation for oscillatory kernels
is presented in this work by developing a curvelet based method.
Here we summarize some of its main features:
\begin{itemize}
	\item The explicitly-sparse representation of the system matrix only consists
		of $O(N \log N)$ nonzero elements.
	\item The computational complexity of the construction of the representation
		with controllable accuracy is log-linear.
	\item S2M, M2M, M2L, L2L, and L2T translation matrices in the directional 
		FMM are used straightforwardly in our curvelet based method.
		As various techniques constructing the low rank approximation of the 
		kernel can be used to compute the translation matrices,
		various variants of our curvelet based method can be easily developed.
	\item It is shown numerically that our method performs well 
		for surface-distributed points with single, 
		double, adjoint, and quadrapole layers, thus it is efficient
		for wave analysis with boundary integral equations as well.
\end{itemize}

Our curvelet based method is constructed as a transform of the directional FMM.
It may also be viewed as the generalization of a wavelet based method to 
high frequency cases, and used as a new wideband fast algorithm.
This work is expected to lay ground to future work related to new fast direct 
solvers and efficient preconditioners for high frequency problems.

\section*{Acknowledgements}

This work is supported by the National Science Foundation of China under Grant 
No. 12101064.
The authors are also greatly appreciated for valuable discussions with 
Prof. Lihua Wen, Jinyou Xiao, Junjie Rong, etc.


\addcontentsline{toc}{chapter}{References}
\bibliography{myref}

\end{document}